\documentclass[11pt]{amsart}
\usepackage{amssymb}
\usepackage{epsfig}
\usepackage{graphicx}

\usepackage{fullpage}
\usepackage[numbers,sort&compress]{natbib}
\usepackage{color}
\usepackage{setspace}
\usepackage[marginal]{footmisc}

\begin{document}
\newtheorem*{pro*}{Proposition~$3.7^{(')}$}
\newtheorem{mthm}{Theorem}
\newtheorem{mcor}{Corollary}
\newtheorem{mpro}{Proposition}
\newtheorem{mfig}{figure}
\newtheorem{mlem}{Lemma}
\newtheorem{mdef}{Definition}
\newtheorem{mrem}{Remark}
\newtheorem{mpic}{Picture}
\newtheorem{rem}{Remark}[section]
\newcommand{\ra}{{\mbox{$\rightarrow$}}}
\newtheorem{Remark}{Remark}[section]
\newtheorem{thm}{Theorem}[section]
\newtheorem{pro}{Proposition}[section]
\newtheorem*{proA}{Proposition A}

\newtheorem{lem}{Lemma}[section]
\newtheorem{defi}{Definition}[section]
\newtheorem{cor}{Corollary}[section]

\title[Classification of solutions]{Classification of solutions to several semi-linear polyharmonic equations and fractional  equations$^*$}
\author[Z. Du, Z. Feng and Y. Li]{Zhuoran Du$^\dag$, Zhenping Feng$^\ddag$ and Yuan Li$^\S$}

\thanks{$^\dag$School of Mathematics, Hunan University, Changsha 410082, China.
{Email: duzr@hnu.edu.cn}. }
\thanks{}

\thanks{$^\ddag$ School of Mathematics, Hunan University, Changsha 410082, China.
{Email: fengzp@hnu.edu.cn}.}
\thanks{}

\thanks{$^\S$ School of Mathematics, Hunan University, Changsha 410082, China.
{Email: liy93@hnu.edu.cn}. }
\thanks{}

\date{\today}

\thanks{}

\date{}

\maketitle

\begin{abstract}
We are concerned with the following semi-linear polyharmonic equation with integral constraint
\begin{align}
\left\{\begin{array}{rl}
&(-\Delta)^pu=u^\gamma_+  ~~ \mbox{ in }{\mathbb{R}^n},\\ \nonumber
&\int_{\mathbb{R}^n}u_+^{\gamma}dx<+\infty,
\end{array}\right.
\end{align}
where $\gamma\in(1,\frac{n+2p}{n-2p})$, $n>2p$, $p\geq2$ and $p\in\mathbb{Z}$. We obtain   that any  nonconstant solution   satisfying certain conditions at infinity is radial symmetric about some point in $\mathbb{R}^{n}$  and monotone decreasing in the radial direction. For the following fractional equation with integral constraint
\begin{equation*}
\left\{\begin{array}{rl}
&(-\Delta)^sv=v^\gamma_+  ~~ \mbox{ in }{\mathbb{R}^n},~~~~\\
&\int_{\mathbb{R}^n}v_+^{\frac{n(\gamma-1)}{2s}}dx<+\infty,~~~~~
\end{array}\right.
\end{equation*}
where $\gamma \in (1, \frac{n+2s}{n-2s})$, $s\in(0,1)$  and $n\geq 2$, we also complete the classification of solutions with certain growth at infinity. In addition, observe that the assumptions of the maximum principle named decay at infinity in \cite{chen} can be weakened slightly. Based on this observation, we classify all positive solutions of two semi-linear fractional  equations  without integral constraint.

\end{abstract}

\noindent
{\it \footnotesize 2020 Mathematics Subject Classification}: {\scriptsize 35B06; 35B08;  35J91.}\\
{\it \footnotesize Key words:  Classification, Polyharmonic equation, Fractional  equation,  Method of moving planes, Radial symmetry.} {\scriptsize }

\section{Introduction}

In this paper, we mainly consider the following two nonlinear equations with integral constraint.
We are devoted to classify their solutions.
One is the polyharmonic equation
\begin{align}\label{01}
\left\{\begin{array}{rl}
&(-\Delta)^pu=u^\gamma_+  ~~ \mbox{ in }{\mathbb{R}^n},\\
&\int_{\mathbb{R}^n}u_+^{\gamma}dx<+\infty,
\end{array}\right.
\end{align}
where $\gamma\in(1,\frac{n+2p}{n-2p})$, $n>2p$, $p\geq2$, $p\in \mathbb{Z}$ and $u_+=\max\{u,0\}$. The other is  the fractional equation
\begin{equation}\label{1}
\left\{\begin{array}{rl}
&(-\Delta)^sv=v^\gamma_+  ~~ \mbox{ in }{\mathbb{R}^n},\\
&\int_{\mathbb{R}^n}v_+^{\frac{n(\gamma-1)}{2s}}dx<+\infty,
\end{array}\right.
\end{equation}
where $s\in(0,1)$, $\gamma \in (1, \frac{n+2s}{n-2s})$ and $n\geq 2$. The fractional Laplacian $(-\Delta)^s$ is defined by
\begin{equation}
\begin{split}
(-\Delta)^s f(x)=C_{s,n} \mbox{P.V.} \int_{{\mathbb{R}}^n}\frac{f(x)-f(y)}{|x-y|^{n+2s}}dy,
\end{split}
\end{equation}
where the constant $C_{s,n}=\frac{2^{2s}\Gamma(\frac{n+2s}{2})}{-\pi^{n/2}\Gamma(-s)}$.
Denote
\begin{equation*}
\mathcal{L}_{2s}({\mathbb{R}}^n):=\left\{f\in L^1_{\mbox{loc}}({\mathbb{R}}^n): \int_{{\mathbb{R}}^n}\frac{|f(x)|}{1+|x|^{n+2s}}dx< \infty\right\}.
\end{equation*}
It is known that $(-\Delta)^sf$ is well defined for $f\in C_{\mbox{loc}}^{2s+\alpha} \cap \mathcal{L}_{2s}({\mathbb{R}}^n)$, where $0<\alpha<1$.
Therefore, throughout the full text, we always assume that  solutions of (\ref{1}) belong to  $ C_{\mbox{loc}}^{2s+\alpha} \cap \mathcal{L}_{2s}({\mathbb{R}}^n)$.

\par  Classification of solutions to a semi-linear elliptic equation is important to understand the considered equation comprehensively and has been extensively studied(see \cite{15}, \cite{Cao}, \cite{CY}, \cite{16}, \cite{chen}, \cite{4}, \cite{Mar}, \cite{Piz}, \cite{WJX}, \cite{XUXING} and the references therein). The method of moving planes ( \cite{15}, \cite{CY}, \cite{21}, \cite{16}, \cite{18}, \cite{23}, \cite{17}), the method of moving spheres (\cite{chenLiz}, \cite{lizhu}) and the sliding methods (\cite{liuz}, \cite{wuchen}) play  important roles in it.

We first recall some results of classification of solutions to several related equations.
For the well-known Yamabe equation
\begin{equation}\label{Y1}
-\Delta u=u^{\frac{n+2}{n-2}} \quad \mbox{in} \quad {\mathbb{R}^n}, ~ n>2,
\end{equation}
Caffarelli, Gidas and Spruck \cite{15} (see also \cite{16} and \cite{18}) classified all its positive solutions by employing the method of moving planes. The natural generalization of (\ref{Y1}) is the  higher order conformally invariant equations
\begin{equation}\label{Y3}
(-\Delta)^p u=u^{\frac{n+2p}{n-2p}}\quad \mbox{in} \quad {\mathbb{R}^n}, ~ n> 2p,~p>0,~p\in \mathbb{Z} .
\end{equation}
Lin \cite{4} proved the classification results for all positive smooth solutions of (\ref{Y3}) in the case of $p=2$. For general integer $p>2$, Wei and Xu \cite{WJX} (see also \cite{XUXING} ) extended the results of \cite{4}.  The classification of positive solutions of (\ref{Y1}) and (\ref{Y3}) is completely solved without adding any assumptions about integrability or decay of solutions. Wang and Ye \cite{12} considered the following problem
\begin{equation} \label{0a0}
-\Delta v=v^{\frac{n}{n-2}}_+  ~~ \mbox{ in }{\mathbb{R}^n},~~n>2,~~\int_{\mathbb{R}^n}v_+^{\frac{n}{n-2}}dx<+\infty,
\end{equation}
and they completed the classification of all nonconstant solutions. Suzuki and Takahashi \cite{1} extended the results of \cite{12} from the Serrin exponent $\frac{n}{n-2}$  to   general subcritical Sobolev exponents. Precisely, they considered the problem
\begin{equation}\label{0a1}
-\Delta v=v^{\gamma}_+  ~~ \mbox{ in }{\mathbb{R}^n},~~n>2,~~\int_{\mathbb{R}^n}v_+^{\frac{n(\gamma-1)}{2}}dx<+\infty,
\end{equation}
where $\gamma\in(1,\frac{n+2}{n-2})$. Chammakhi, Harrabi and Selmi \cite{RHA} classified  all  sign-changing  solutions of
\begin{equation}\label{0a2}
\Delta^2 v=v^{\gamma}_+  ~~ \mbox{ in }{\mathbb{R}^n},~~1<\gamma\leq \frac{n}{n-4},~~n>4,~~
\int_{\mathbb{R}^n}v_+^{\gamma}dx<+\infty.
\end{equation}

Classification of solutions to fractional equations is also drawing many researchers  recently. By developing  a direct method of moving planes for  fractional Laplace equations, Chen, Li and Li in \cite{chen} considered the equation
$$ (-\Delta)^{s}u=u^{\frac{n+2s}{n-2s}}~~~ ~\mbox{in} ~~{\mathbb{R}^n}, ~~~0<s<1.$$
They obtained that its any positive solution must be radially symmetric about some point in ${\mathbb{R}^n}$ and monotone decreasing in radial direction, which is a generalization of the result in \cite{15}. Dai and Qin in \cite{Qin} derived the classification of nonnegative classical solutions to the following problem
$$ (-\Delta)^{\frac{3}{2}}u=u^{\frac{n+3}{n-3}}~~~ ~\mbox{in} ~~{\mathbb{R}^n},~~n>3,~~\int_{\mathbb{R}^n}\frac{u^{\frac{n+3}{n-3}}}{|x|^{n-3}}dx< \infty.$$
Recently, Cao, Dai and Qin \cite{Cao} completed the classification of nonnegative classical solutions to the  higher-order fractional Laplace equation
$$ (-\Delta)^{m+\frac{\alpha}{2}}u=u^{\frac{n+2m+\alpha}{n-2m-\alpha}}~~~ ~\mbox{in} ~~{\mathbb{R}^n},~~n>2m+\alpha,~~m\geq 1,~0<\alpha<2,$$
which  improve the classification results in \cite{Qin} by removing the integrability assumption.

Inspired by the above works, we focus on the classification of bounded energy solutions and want to extend the results of (\ref{0a0})-(\ref{0a2}) in two aspects. Precisely, we will extend the results of (\ref{0a0}) and (\ref{0a2}) to polyharmonic equations with any subcritical Sobolev exponents with respect to the space $H^p(\mathbb{R}^n)$, and the results of (\ref{0a1}) to  fractional equations.

Our results are as follows.

\begin{thm}\label{T1}
Let $\gamma\in(1,\frac{n+2p}{n-2p})$, $n>2p$ and $p\geq2$. Suppose that $u\in C^{2p}({\mathbb{R}^n})$ is a nonconstant solution of (\ref{01}) and satisfies $u(x)=o(|x|^2)$ at infinity, then $u$ is symmetric about some point $x_0\in{\mathbb{R}^n}$ and $\frac{\partial u}{\partial r}<0$, where $r=|x-x_0|$.
\end{thm}

Note that Theorem \ref{T1} extends the results of \cite{RHA} to all subcritical Sobolev exponent $\gamma \in (1,\frac{n+4}{n-4})$.

From the nonexistence of positive solutions in subcritical Sobolev exponent case of \cite{WJX}, we know that (\ref{01}) does not possess positive solutions.
It is easy to verify that $u(x)=-\sum_{i=1}^n a_ix_i^2$ with all $a_i\geq 0$   are solutions of (\ref{01}). Plainly the growth assumption $u(x)=\mbox{o}(|x|^2)$ at infinity  rules out these solutions.
In fact, the growth assumption $u(x)=\mbox{o}(|x|^2)$  can rule out all negative nonconstant solutions.
 Suppose (\ref{01}) has a negative nonconstant solution $u$. Then the equation in (\ref{01}) tells us that $(-\Delta)^{p-1}u$ is harmonic. By Lemma \ref{a1}, we have $(-\Delta)^{p-1}u\geq 0$ in ${\mathbb{R}^n}$. Liouville's Theorem  implies that
$$ (-\Delta)^{p-1}u\equiv c\geq 0.$$
From this and $\lim\limits_{|x|\rightarrow \infty}(-\Delta)^{p-1}u=0$ (see (\ref{y1})), we derive $c=0$. Since $(-\Delta)^{j}u\geq 0$ and $\lim\limits_{|x|\rightarrow \infty}(-\Delta)^{j}u=0$ hold for any $j=1,\ldots,p-1$, after repeating $p-1$ times of the above argument, we have $-\Delta u=0$ in ${\mathbb{R}^n}$. Due to $u$ is a negative solution so we deduce that $u$ must be a constant, which is a contradiction. Hence any nonconstant solution of (\ref{01}) satisfying the growth $\mbox{o}(|x|^2)$ at infinity must be sign-changing.

Under the integral constraint, for any solution of (\ref{01}) we can prove
that the decay assumption $\Delta u(x)\rightarrow 0$ as $|x|\rightarrow \infty$ is equivalent to the growth assumption $u(x)=\mbox{o}(|x|^2)$ at infinity. Therefore we have the following result.

\begin{thm}\label{T111}
Let $\gamma\in(1,\frac{n+2p}{n-2p})$, $n>2p$ and $p\geq2$. Suppose that $u\in C^{2p}({\mathbb{R}^n})$ is a nonconstant solution of (\ref{01}) and satisfies $\Delta u(x)\rightarrow 0$ as $|x|\rightarrow \infty$, then $u$ is symmetric about some point $x_0\in{\mathbb{R}^n}$ and $\frac{\partial u}{\partial r}<0$, where $r=|x-x_0|$.
\end{thm}

The forthcoming paper \cite{Du} deals with the classification of all nonconstant solutions to the problem (\ref{01}) with higher-order fractional Laplacians $(-\Delta)^{p+s}$, where  $n> 2(p+s)$, $0<s<1$ and $p\geq 1$ is an integer.

Our results on the fractional problem (\ref{1}) are as follows.

\begin{thm}\label{F1}
Assume that $\gamma \in (1, \frac{n+2s}{n-2s})$, $n\geq 2$ and $\tau<2s$. Let $v$ be a nonconstant  solution of (\ref{1}) satisfying $|v(x)|= O(|x|^\tau)$ at infinity. Then $v$ is symmetric about some point $x_0\in{\mathbb{R}^n}$ and $\frac{\partial v}{\partial r}<0$.
\end{thm}

Here $u(x)=\mbox{o}(|x|^2)$ at infinity means that for any $\varepsilon>0$, there exists that $R>0$ such that $|u(x)|\leq \varepsilon |x|^2$ for $|x|\geq R$. And $|v(x)|= \mbox{O}(|x|^\tau)$ at infinity means that there exists $R>0$ and $C_0>0$ such that
$\frac{C_0}{2}|x|^\tau \leq |v(x)| \leq 2 C_0 |x|^\tau \mbox{ for any}~ |x|\geq R.$

From the nonexistence of positive solutions in the subcritical exponent case of \cite{chen} and the Liouville theorem of $s$-harmonic function of \cite{8}, it is easy to verify that any nonconstant solution of (\ref{1}) is necessarily a sign-changing solution.

Note that (\ref{1}) is invariant under the scaling $\mu^{\frac{2s}{\gamma-1}}v(\mu x), \mu>0$, which will be used in Lemma 3.2 (Section 3).
If the integral assumption $\int_{\mathbb{R}^n}v_+^{\frac{n(\gamma-1)}{2s}}dx<+\infty$ in (\ref{1}) is revised into $\int_{\mathbb{R}^n}v_+^{\gamma}dx<+\infty$, which seems more natural for the equation in (\ref{1}),  then the above scaling invariance shows that only value of $\gamma$ is $\frac{n}{n-2s}$.
So the  integral assumption $\int_{\mathbb{R}^n}v_+^{\gamma}dx<+\infty$ only corresponds to the particular case $\gamma=\frac{n}{n-2s}$.
Note that the solution of the equation $\frac{n(\gamma-1)}{2s}=\gamma$ is $\gamma=\frac{n}{n-2s}$ exactly.
Hence, for general $\gamma\in (1, \frac{n+2s}{n-2s})$,   the integral assumption  $\int_{\mathbb{R}^n}v_+^{\frac{n(\gamma-1)}{2s}}dx<+\infty$ seems to be an  appropriate candidate.

Meantime, one may note that an integral assumption is included in (\ref{01}), and is not needed in the classification of solutions to (\ref{Y3}) with a pure power critical nonlinearity.
We believe that the possible reason of this difference is that the equation  (\ref{01}) is much closer to  the following Liouville equation (see \cite{WJX}) than (\ref{Y3}) is
\begin{align*}\label{0abc11}
\left\{\begin{array}{rl}
&(-\Delta)^p u=(2p-1)!e^{2pu}  ~~ \mbox{ in }{\mathbb{R}^{2p}},\\
&\int_{\mathbb{R}^{2p}}e^{2pu}dx<+\infty.
\end{array}\right.
\end{align*}
Similar view (the case $p=1$) was pointed out by Wang and Ye in \cite{12}.

 The authors in \cite{chen} obtained  a key ingredients (Decay at infinity, see Proposition \ref{FL7} in Section 3 of this paper), which is crucial to establishing the method of moving planes for  fractional Laplace equations in unbounded domains, through the integral defining of the fractional laplacian.
In the light of this, they showed that any positive solution of the nonlinear Schr\"{o}dinger equation with fractional diffusion
  \begin{equation} \label{fse}
 (-\Delta)^su+u=u^\nu,~ 1<\nu<\infty~~x\in \mathbb{R}^n,
 \end{equation}
 satisfying $\lim\limits_{|x|\rightarrow\infty}u(x)=a<\left(\frac{1}{\nu}\right)^\frac{1}{\nu-1}$, must be radially symmetric about some point in $\mathbb{R}^n$
 and decreasing in the radial direction \cite{chen}.

 Observe that the result of Proposition \ref{FL7} still holds true after weakening its assumptions, and we establish Proposition $3.7^{(')}$(see Section 3). Further, we  improve the above result of (\ref{fse}) in \cite{chen}
 as follows.
\begin{thm}\label{F2}
Suppose $C^1$ function $g$ satisfies
\begin{equation}\label{1.6}
|g'(u)|\leq C |u|^\mu
\end{equation}
for some $\mu>0, C>0.$ Assume that $u\in \mathcal{L}_{2s}\cap C_{\mbox{loc}}^{2s+\alpha}$ is a positive solution of
\begin{equation*}
(-\triangle)^su+u=g(u),~~x\in\mathbb{R}^n
\end{equation*}
satisfying \begin{equation}\label{1.7}
u(x)< \left(\frac{\tilde{C}|x|^{-2s}+1}{C}\right)^{\frac{1}{\mu}}~ \mbox{ for } |x| \mbox{ large enough },
\end{equation}
where $\tilde{C}=\frac{2w_n C_{s,n}}{4^{n+2s}}$, then $u$ must be radially symmetric  about some point in ${\mathbb{R}^n}$ and monotone decreasing in the radial direction.
\end{thm}
\begin{rem}\label{FR1}
\rm Equation (\ref{fse}) corresponds to the particular case $g(u)=u^\nu$ of Theorem \ref{F2}. Note that for this $g$ we have $\mu=\nu-1>0, C=\nu$
and
$$
\lim_{|x|\rightarrow\infty} \left(\frac{\tilde{C}|x|^{-2s}+1}{C}\right)^{\frac{1}{\mu}}=\left(\frac{1}{\nu}\right)^{\frac{1}{\nu-1}}.
$$
 Obviously, the assumption (\ref{1.7}) is satisfied if $\lim\limits_{|x|\rightarrow \infty}u(x)=a < (\frac{1}{\nu})^{\frac{1}{\nu-1}}$.
\end{rem}

Applying Proposition $3.7^{(')}$, we can also classify all positive solutions of the following semilinear fractional Laplace equation  without integral constraint.
\begin{thm}\label{F3}
Let $\mu>0$, $C>0$,  $n\geq2$ and $h$ satisfies
\begin{equation}\label{1.9}
\left|\frac{h(u)-h(v)}{u-v}\right|\leq C|u+v|^\mu \mbox{ for all } u\neq v.
\end{equation}
Assume that $u\in \mathcal{L}_{2s}\cap C_{\mbox{loc}}^{2s+\alpha}$ is a positive  solution of the problem
\begin{equation}\label{1.8}
\begin{cases}
(-\Delta)^s{u}=h(u) &\mbox{ in }{{\mathbb{R}}^n},\\
\lim\limits_{|x|\rightarrow \infty}u(x)=0. &\\
\end{cases}
\end{equation}
If $u(x)=O(\frac{1}{|x|^m})$ at infinity for some $m$ satisfying  $m >\frac{2s}{\mu} $, then $u$ is radially symmetric about some point in ${\mathbb{R}}^n$ and monotone decreasing in the radial direction.
\end{thm}

Theorem \ref{F3} is an improvement  of Theorem 1.2 in \cite{14}, where similar result is established  for  $m>\max \{\frac{2s}{\mu},\frac{n}{\mu+2}\}$.
Similar results of Theorems \ref{F2} can be found in \cite{13}.

\par This paper is organized as follows. In Section 2, we shall complete the proof of Theorem \ref{T1} and Theorem \ref{T111}.  In Section 3, we will first give some preliminary results and lemmas, and then complete the classification of (\ref{1}). In the last section, we will use the moving plane method of fractional Laplacian to finish the proofs of Theorems \ref{F2} and \ref{F3}.

\section{Radial symmetry of solutions to polyharmonic equations}
In this section, we will complete the proof of Theorem \ref{T1} and Theorem \ref{T111}. In order to prove Theorems \ref{T1} and \ref{T111}, a key ingredient is to prove that solution $u$ is bounded from above.
For this purpose, we need to establish several technique lemmas.

Let $u\in C^{2p}({\mathbb{R}^n})$ be a solution of (\ref{01}) and set
\begin{align}\label{ppffaa1}
\zeta(x)=-C_{n,p}\int_{\mathbb{R}^n}\frac{u_+^\gamma(y)}{|x-y|^{n-2p}}dy,
\end{align}
where $C_{n,p}=\frac{2\Gamma(\frac{n}{2}-p)}{ne_n4^p\Gamma(\frac{n}{2})(p-\frac{n}{2})!(p-1)!}$ and $e_n=\frac{\pi^{\frac{n}{2}}}{\Gamma(1+\frac{n}{2})}$.
Then $\zeta$ is well-defined and $\zeta\in C^{2p}({\mathbb{R}^n})$.

\subsection{Growth assumption $u(x)=\mbox{o}(|x|^2)$ at infinity}$\\$

We first establish the super polyharmonic property of equation (\ref{01}).
\begin{lem}\label{a1}
Suppose that $\gamma\in(1,\frac{n+2p}{n-2p})$, $n>2p$, $p\geq 2$ and $u$ is a solution of equation (\ref{01}) satisfying $u(x)=o(|x|^2)$ at infinity, then we have $$(-\Delta)^i u\geq 0, ~~~i=1,2,...,p-1.$$
\end{lem}
{\bf Proof } Denote $w_i=(-\Delta)^iu$, $i=1,2,...,p$. Firstly, we prove that $w_{p-1}\geq 0$ by contradiction.
\par Suppose there exists $x_0 \in {\mathbb{R}^n}$ such that $w_{p-1}(x_0)< 0$. Without loss of generality, we assume that $x_0=0$.
For continuous function $f$, define
$$ \bar{f}(r):=\frac{1}{|\partial B_r(0)|}\int_{\partial B_r(0)} fd\sigma, ~~~\mbox{ for all }~ r>0,$$
and $ \bar{f}(0):=\lim\limits_{r\rightarrow0}\bar{f}(r)$, then $\bar{f}(0)=f(0)$. By the well-known property
$\Delta \bar{f}(r)=\overline{\Delta f}(r)$, we have
\begin{align}
\begin{split}
&\Delta  \bar{u}+\bar{w}_1=0,\nonumber\\
&\Delta  \bar{w}_1+\bar{w}_2=0,\nonumber\\
&...\nonumber\\
&\Delta \bar{w}_{p-1} +\bar{u}_+^\gamma\leq 0.
\end{split}
\end{align}
where we used Jensen's inequality in the last formula. From the last inequality, we can see that $ \bar{w}'_{p-1}(r)\leq0$ for all $r\in (0,+\infty)$. Combining  this with $ \bar{w}_{p-1}(0)=w_{p-1}(0)<0$, we obtain that
\begin{equation}\label{X1}
\bar{w}_{p-1}(r)\leq \bar{w}_{p-1}(0)<0~  \mbox{ for all } r.
\end{equation}
It is easy to see that
$$ \Delta \bar{w}_{p-2}(r)=-\bar{w}_{p-1}(r)\geq -\bar{w}_{p-1}(0),$$
it follows that
$$ \bar{w}'_{p-2}(r)\geq \frac{-\bar{w}_{p-1}(0)}{n}r=:cr,  $$
where $c$ is a positive constant. By simple calculation we deduce that there exist $c_2,r_1>0$ such that
\begin{equation}
\bar{w}_{p-2}(r)\geq \frac{c}{2}r^2+\bar{w}_{p-2}(0)\geq \frac{c}{4}r^2=:c_2 r^2~  \mbox{ for } r\geq r_1.
\end{equation}
Same arguments give that there exist $c_3>0$ and $r_2\geq r_1$ such that
\begin{equation}
\bar{w}_{p-3}(r)\leq -c_3 r^4~ \mbox{ for } r\geq r_2,
\end{equation}
and there exist $c_i>0$ and $r_{i-1}\geq r_{i-2}$ $(i=1,...,p)$ such that
\begin{equation}\label{20220721e}
(-1)^i\bar{w}_{p-i}(r)\geq c_i r^{2(i-1)}~ \mbox{ for } r\geq r_{i-1}.
\end{equation}
In particular,
$$ (-1)^p\bar{u}(r)\geq c_p r^{2(p-1)}~ \mbox{ for}~ r\geq r_{p-1}.  $$
Since $p\geq 2$, we conclude $|\bar{u}(r)|\geq c_p r^2$ for $r\geq r_{p-1}$, which  contradicts  with $u(x)=\mbox{o}(|x|^2)$. Hence $w_{p-1}\geq 0$.

The results $w_{p-i} \geq 0$ for all $i=1,2,...,p-1$ now can be proved by mathematical induction.
For $i=1$, the above argument have shown  $w_{p-1}\geq 0$. We assume that $w_{p-i} \geq 0$ for all $i\leq k\leq p-2$.
Suppose $w_{p-(k+1)}<0$ for some $x_0$, and we may assume that $x_0=0$. By assumption $w_{p-k}\geq 0$, applying the above system again, we obtain
$\Delta \bar{w}_{p-(k+1)}<0$. Integrate it to get that $\bar{w}_{p-(k+1)}\leq\bar{w}_{p-(k+1)}(0)<0$.
Similar arguments as above show that $\bar{w}_{p-(k+2)}\geq c_2r^2$ for some $c_2$.
Repeatedly integrating the system we obtain that for large $r$
$$
(-1)^j\bar{w}_{p-(k+j)}(r)\geq c_j r^{2(j-1)}, ~~~k+j\leq p.
$$
In particular, choose $k+j=p$, then $j\geq 2$, since $k\leq p-2$.
Thus the above inequality gives that $|\bar{u}(r)|\geq cr^2$ at infinity for some $c>0$, which contradicts  with $u(x)=\mbox{o}(|x|^2)$ again.
$\hfill\square$
\begin{lem}\label{0723abc}
Suppose that $\gamma\in(1,\frac{n+2p}{n-2p})$, $n>2p$, $p\geq 2$ and $u$ is a solution of equation (\ref{01}) satisfying $u(x)=o(|x|^2)$ at infinity. If $\zeta(x)$ is given by (\ref{ppffaa1}), then $(-\Delta)^i(u+\zeta)=0$, $i=1,2,...,p-1$.
\end{lem}
{\bf Proof } As in \cite{220728} we introduce  the following boundary value problem
\begin{align}\label{f2207271}
\begin{cases}
& (-\Delta)^{p-1} \phi =\delta(x)\ \ \ \ \ \ \ \ \ \ \ \ \ \  \ \mbox{ in }  B_r(0),\\
&  \phi=\Delta \phi=\cdot \cdot \cdot=\Delta ^{p-2} \phi =0\ \ \mbox{ on } \partial B_r(0),
\end{cases}
\end{align}
where $\delta (x)$ is the Dirac Delta function. By the maximum principle, for the unique solution $\phi$ of (\ref{f2207271}) we can easily obtain that
\begin{align}\label{220727y}
|\phi(x)|\leq \frac{C}{|x|^{n-2(p-1)}},
\end{align}
and for any $j=0,1,...,p-2$,
\begin{equation}\label{f2207272}
\frac{\partial}{\partial \nu} [(-\Delta)^j \phi]\leq 0~~\mbox{ on } \partial B_r(0).
\end{equation}
One can verify that
\begin{align}\label{220728a}
\left|\frac{\partial}{\partial\nu}[(-\Delta)^j \phi]\right|\leq \frac{C}{r^{n-2p+3+2j}}~~\mbox{ on } \partial B_r(0).
\end{align}

Denote $v=-\Delta u$, then
\begin{align}\label{f2207273}
 (-\Delta)^{p-1} v(x)=u_+^\gamma(x)~~~\mbox{ in } \mathbb{R}^n.
\end{align}
Multiplying  both side of the equation (\ref{f2207273}) by $\phi$ and integrating by parts several times on $B_r(0)$, we arrive at
\begin{align}\label{f2207274}
\begin{split}
\int_{B_r(0)}u_+^{\gamma}(x)\phi(x)dx&=\int_{B_r(0)}(-\Delta)^{p-1}v(x)\phi(x)dx\\
&=v(0)+\sum\limits_{j=0}^{p-2}\int_{ \partial B_r(0)}(-\Delta)^{j}v(x)\frac{\partial}{\partial \nu}[(-\Delta)^{p-2-j}\phi(x)]d\sigma.
\end{split}
\end{align}
Applying Lemma \ref{a1} and (\ref{f2207272}), we have
\begin{align}\label{f2207276}
\int_{B_r(0)}u_+^{\gamma}(x)\phi(x)dx \leq   v(0).
\end{align}
Obviously
\begin{align}\label{220728e}
\phi(x)\rightarrow \frac{C_{n,p-1}}{|x|^{n-2(p-1)}}~~~\mbox{ as } r\rightarrow \infty,
\end{align}
where $C_{n,p-1}=C_{n,p}(2p-2)(n-2p)$. This and (\ref{f2207276}) yield that there is a constant $c(n)$ such that
\begin{align}\label{220727e}
\int_{\mathbb{R}^n}\frac{u_+^\gamma(x)}{|x|^{n-2(p-1)}}dx\leq c(n)v(0)<+\infty.
\end{align}
Hence there exists $r_m\rightarrow \infty$ such that
$$\int_{\partial B_{r_m}(0)}\frac{u_+^\gamma(x)}{|x|^{n-2(p-1)-1}}d\sigma \rightarrow 0. $$

Applying the above arguments to the following equation instead of (\ref{f2207273})
$$ (-\Delta)^{i-1} v=(-\Delta)^{i}u$$
for each $i=2,..., p-1$, we have
$$ \int_{\mathbb{R}^n}\frac{(-\Delta)^{i}u}{|x|^{n-2(i-1)}}dx<+\infty.$$
Hence there exists $r_m\rightarrow \infty$ such that
\begin{align}\label{220727i}
\sum_{i=1}^{p-2}\int_{\partial B_{r_m}(0)} \frac{(-\Delta)^{i} v}{|x|^{n-2i-1}}d\sigma \rightarrow 0.
\end{align}

Next we claim that
\begin{align}\label{220727p}
\int_{\partial B_{r_m}(0)} \frac{ v}{|x|^{n-1}}d\sigma \rightarrow 0, ~~~\mbox{ as } r_m\rightarrow \infty.
\end{align}
To prove (\ref{220727p}), it is enough to show that $\lim\limits_{r\rightarrow \infty} \bar{v}(r)=0$.
Otherwise, there exist $\varepsilon_0>0$ such that for any positive integer $k$, we have $r_k>k$ satisfying $|\bar{v}(r_k)|\geq \varepsilon_0$.
 Note that $\Delta\bar{v}(r)\leq0$ for any $r\geq0$, from the maximum principle we have that $\bar{v}(r)\geq \varepsilon_0$ for any $r\leq r_k$. Hence $\Delta \bar{u}(r)\leq -\varepsilon_0$, which yields
$$ \bar{u}(r)\leq \bar{u}(0)-\frac{\varepsilon_0}{2n}r^2, ~~r\leq r_k.   $$
This contradicts with the growth assumption $u(x)=\mbox{o}(|x|^2)$ at infinity. Hence the claim holds.

By virtue of (\ref{220728a}), (\ref{220727i}) and (\ref{220727p}), there exist a sequence $r_m\rightarrow \infty$ such that the boundary integrals  on $\partial B_{r_m}(0)$ in (\ref{f2207274}) approach 0 as $r_m\rightarrow \infty$. From (\ref{220727y}), (\ref{220728e}) and (\ref{220727e}), by using the Lebesgue Dominated Convergence Theorem to the left hand side of (\ref{f2207274}), and taking limit along the sequence $r_m\rightarrow \infty$, we conclude that
$$ -\Delta u(0)=C_{n,p-1}\int_{\mathbb{R}^n}\frac{u_+^\gamma(y)}{|y|^{n-2{p}+2}}dy .   $$
By a translation, we derive for any $x\in \mathbb{R}^n$,
$$ -\Delta u(x)=C_{n,p-1}\int_{\mathbb{R}^n}\frac{u_+^\gamma(y)}{|x-y|^{n-2{p}+2}}dy , $$
that is
$$ -\Delta(u+\zeta)=0,~~~~ \mbox{ in } \mathbb{R}^n.  $$

For $i=2,...,p-1$, we  denote $v=(-\Delta)^i u$ and consider the following boundary value problem instead of (\ref{f2207271})
\begin{align*}
\begin{cases}
& (-\Delta)^{p-i} \phi =\delta(x)\ \ \ \ \ \ \ \ \ \ \ \ \ \  \ \mbox{ in }  B_r(0),\\
&  \phi=\Delta \phi=\cdot \cdot \cdot=\Delta ^{p-i-1} \phi =0\ \ \mbox{ on } \partial B_r(0).
\end{cases}
\end{align*}
Repeating  the above argument we  see that the results $(-\Delta)^i(u+\zeta)=0$, $i=2,...,p-1$ follow. $\hfill\square$

Borrowing the idea of \cite{Guo}, we establish Lemma \ref{z2} and Lemma \ref{0724}.

\begin{lem}\label{z2}
Let $\gamma\in(1,\frac{n+4}{n-4})$, $p=2$ and $n>4$. Suppose that $u$ is a solution of (\ref{01}) satisfying $u(x)=o(|x|^2)$, then  $\Delta u(x)\rightarrow 0$ as $|x|\rightarrow \infty$.
\end{lem}

{\bf Proof } To prove this lemma, noting  that $\limsup\limits_{|x|\rightarrow \infty} \Delta u(x)\leq 0$ as given by  Lemma \ref{a1}, we only need to show $\liminf\limits_{|x|\rightarrow \infty}\Delta u(x)\geq 0$.

Suppose that the assertion fails, then there exist $\varepsilon>0$ and a sequence $\{x_k\}\subset \mathbb{R}^n$ with $|x_k|\rightarrow \infty$ as $k\rightarrow \infty$ such that
$$ \Delta u(x_k)\leq -\varepsilon ~~~~\mbox{ for } k\geq 1.$$
Let $v_k(y)=u(x)$, where $y=x-x_k$, we have
$$ \Delta_y^2v_k=(v_k)_+^\gamma~~\mbox{ and } ~~\Delta_yv_k(0)=\Delta_xu(x_k)\leq -\varepsilon.$$
Set
$$ w_k(y)=\frac{\Delta v_k(y)}{\Delta v_k(0)}.$$
It is easy to see that
$$ w_k(0)=1  ~~\mbox{ and } ~~ \Delta \bar{w}_k(r)=\frac{\overline{(v_k)_+^\gamma}}{\Delta v_k(0)}.$$
Simple calculation shows that
\begin{equation}\label{r3}
r^{n-1}\bar{w}'_k(r)=\frac{1}{\Delta v_k(0)\omega_n}\int_{B_r(0)}(v_k)_+^\gamma(y) dy<0,
\end{equation}
where $\omega_n$ is the surface area of the $n$-dimensional unit sphere. For any fixed $0<R\leq \frac{1}{2}|x_k|$, we have
$$ \int_{B_R(0)}(v_k)_+^\gamma(y)dy\leq \int_{B_{\frac{|x_k|}{2}}(x_k)}(u)_+^\gamma(y)dy\rightarrow 0 \mbox{ as } k\rightarrow \infty.$$
Therefore, it follows from (\ref{r3}) that
$$ \bar{w}'_k(r)\rightrightarrows 0 \mbox{ for } r\in (0,R] \mbox{ as } k \rightarrow \infty,$$
which implies
\begin{equation}\label{r5}
\bar{w}_k(r)\rightrightarrows 1 \mbox{ for } r\in [0,R] \mbox{ as } k \rightarrow \infty.
\end{equation}
Here the sign $\rightrightarrows$ means uniform convergence. On the other hand, for $r\in[R,\frac{1}{2}|x_k|]$, we have
$$ -\bar{w}'_k(r)\leq \frac{r^{1-n}}{|\Delta v_k(0)|\omega_n}\int_{{B_{\frac{|x_k|}{2}}}(0)}(v_k)_+^\gamma dy.$$
Integrating on $[R,r]$ yields that for any $r\in [R,\frac{1}{2}|x_k|]$,
\begin{equation*}\label{r6}
0\leq \bar{w}_k(R)-\bar{w}_k(r)\leq \frac{R^{2-n}}{(n-2)|\Delta v_k(0)|\omega_n}\int_{{B_{\frac{|x_k|}{2}}}(0)}(v_k)_+^\gamma dy\rightrightarrows 0~\mbox{ as } k\rightarrow\infty .
\end{equation*}
Hence for any $r\in [0,\frac{1}{2}|x_k|]$, we derive
\begin{equation}\label{r4}
\bar{w}_k(r)\rightrightarrows 1~~\mbox{ as } k\rightarrow\infty .
\end{equation}
From (\ref{r4}), for $k$ large enough and $r\in [0,\frac{1}{2}|x_k|]$, we derive
$$ \Delta \bar{v}_k(r)=\bar{w}_k(r)\Delta v_k(0)\leq \frac{1}{2}\Delta v_k(0)\leq-\frac{1}{2}\varepsilon .$$
Direct integration shows that
$$ \bar{v}_k(\frac{1}{2}|x_k|)-\bar{v}_k(0)\leq -\frac{\varepsilon}{4n}(\frac{1}{2}|x_k|)^2,$$
which yields
$$ \bar{v}_k(\frac{1}{2}|x_k|)\leq -(\frac{\varepsilon}{4n}+o(1))(\frac{1}{2}|x_k|)^2.$$
This contradicts with the assumption $u(x)=\mbox{o}(|x|^2)$ at infinity. Hence $\liminf\limits_{|x|\rightarrow \infty}\Delta u(x)\geq 0$.
$\hfill\square$

Compared with Lemma \ref{z2}, for the cases $p>2$ we only obtain the weaker conclusion that $\Delta u$ is bounded.

\begin{lem}\label{0724}
Suppose that $\gamma\in(1,\frac{n+2p}{n-2p})$, $n>2p$, $p>2$ and $u$ is a solution of equation (\ref{01}) satisfying $u(x)=o(|x|^2)$ at infinity. Then $\Delta u\in L^\infty(\mathbb{R}^n)$.
\end{lem}
{\bf Proof } We will also carry out our proof  via contradiction arguments. Suppose that this lemma does not hold, then there exists a sequence $\{x_k\}\subset \mathbb{R}^n$ with $|x_k|\rightarrow \infty$ as $k\rightarrow \infty$ such that
$$ -\Delta u(x_k)> k.$$
Let $v_k(y)$, $w_k(y)$ and $\bar{w}_k(r)$ be defined as in Lemma \ref{z2}, where $y=x-x_k$ and $r\in[0,{\frac{1}{2}}|x_k|]$. Then we have
$$ \Delta_y^2v_k(y)=\Delta_x^2u(x)~~\mbox{ and } ~~-\Delta_yv_k(0)=-\Delta_xu(x_k)>k.$$
Obviously,
$$ w_k(0)=1  ~~\mbox{ and } ~~ \Delta \bar{w}_k(r)=\frac{{\frac{1}{|\partial {B_r(0)}|}\int_{\partial {B_r(0)}}\Delta^2 u(y+x_k)d\sigma}}{\Delta v_k(0)}.$$
Simple calculation gives
\begin{equation}\label{220727a}
r^{n-1}\bar{w}'_k(r)=\frac{1}{\Delta v_k(0)w_n}\int_{B_r(0)}\Delta^2 u(y+x_k) dy.
\end{equation}
Note that for any $s>0$, $x_0\in \mathbb{R}^n$ and $R>0$, we have
$$ \int_{B_R(x_0)}\frac{dx}{|x|^{n-s}}\leq C(n)R^s,$$
where the constant $C(n)$ depends only on $n$.
Combining this, Lemma \ref{0723abc} and the assumption $p>2$,
for any fixed $0<R<\frac{1}{2}|x_k|$, we deduce
\begin{align}\label{220727a1}
\begin{split}
\int_{B_R(0)}\Delta^2 u(y+x_k)dy
&=\tilde{C}_{n,p}\int_{B_R(x_k)}\int_{\mathbb{R}^n}\frac{u_+^\gamma(x)}{|x-y|^{n-2p+4}}dxdy\\
&=\tilde{C}_{n,p}\int_{\mathbb{R}^n}u_+^\gamma(x)\int_{B_R(x_k+x)}\frac{1}{|y|^{n-2p+4}}dydx\\
&\leq CR^{2p-4}\int_{\mathbb{R}^n}u_+^\gamma(x)dx\leq CR^{2p-4},
\end{split}
\end{align}
where $\tilde{C}_{n,p}=C_{n,p}(n-2p)(2p-2)(2p-4)(n+2-2p)$. The last inequality holds due to the integral constraint in (\ref{01}).
Combining (\ref{220727a})-(\ref{220727a1}) and the fact $\frac{1}{\Delta v_k(0)}\rightarrow 0$ as $k\rightarrow \infty$, we obtain
$$ \bar{w}'_k(r)\rightrightarrows0, ~~~\mbox{ for } r\in(0,R] \mbox{ as } k\rightarrow \infty,$$
which implies
$$ \bar{w}_k(r)\rightrightarrows1,~~~\mbox{ for } r\in(0,R] \mbox{ as } k\rightarrow \infty.$$

For $r\in [R, \frac{1}{2}|x_k|]$, using the same proof as ({\ref{220727a1}}) we have
\begin{align}
\begin{split}
-\bar{w}'_k(r)=\frac{r^{1-n}\int_{B_r(x_k)}\Delta^2u(x)dx}{|\Delta v_k(0)| w_n}\leq \frac{Cr^{1-n}r^{2p-4}}{|\Delta v_k(0)| w_n}=\frac{Cr^{2p-n-3}}{|\Delta v_k(0)| w_n}.
\end{split}
\end{align}
Integrating on $[R,r]$ yields that for any $r\in [R,\frac{1}{2}|x_k|]$,
\begin{equation*}\label{r6}
0\leq \bar{w}_k(R)-\bar{w}_k(r)\leq \frac{C R^{2p-n-2}}{(n+2-2p)|\Delta v_k(0)|w_n}\rightrightarrows 0~\mbox{ as } k\rightarrow\infty .
\end{equation*}
Hence for any $r\in [0,\frac{1}{2}|x_k|]$, we have
\begin{equation}\label{r4q}
\bar{w}_k(r)\rightrightarrows 1~~\mbox{ as } k\rightarrow\infty .
\end{equation}
From (\ref{r4q}), we derive for $k$ large enough and $r\in [0,\frac{1}{2}|x_k|]$,
$$ \Delta \bar{v}_k(r)\leq \frac{1}{2}\Delta v_k(0)< -\frac{1}{2}k .$$
The rest argument that lead to a contradiction with the growth assumption at infinity is similar to that of Lemma \ref{z2} and we omit it.
$\hfill\square$

From Lemmas \ref{z2} and  \ref{0724}, we have $\Delta u\in L^\infty(\mathbb{R}^n)$. Now we can prove that $u$ is bounded from above.

\begin{lem}\label{a2}
Assume $\gamma\in(1,\frac{n+2p}{n-2p})$, $n>2p$, $p\geq2$ and $u$ is a solution of (\ref{01}) satisfying $u(x)=o(|x|^2)$ at infinity, then there exists a constant $C>0$ such that $\sup\limits_{\mathbb{R}^n}u \leq C$.
\end{lem}
{\bf Proof} From Lemmas \ref{z2} and \ref{0724}, we have that there exists $A>0$ such that $|\Delta u| \leq A$ in ${\mathbb{R}^n}.$
Lemma \ref{a1} implies that $\Delta u\leq 0$. Denote $h(x)=-\Delta u(x)$. Given
$x_0\in {\mathbb{R}^n}$, let $u_1$ be the solution of
\begin{align}\label{03}
\left\{\begin{array}{rl}
(-\Delta)v=h ~~~&\mbox{ in } B_1(x_0),\\
v=0~~~~~~~~~~ &\mbox{ on } \partial B_1(x_0).
\end{array}\right.
\end{align}
It follows from the elliptic theory that $|u_1|\leq C$, where $C>0$ independent of $x_0$. Set $u_2=u-u_1$,
then $(u_2)_+\leq u_++|u_1|$. Since $|u_1|\leq C$ in $B_1(x_0)$ and $\int_{\mathbb{R}^n}u_+^{\gamma}(x)dx<+\infty$, we derive
$$\int_{B_1(x_0)}(u_2)_+^{\gamma}(x)\leq C.$$
Note that $\Delta u_2=0$ in $B_1(x_0)$.
For the subharmonic function $(u_2)_+$, we have
$$\|(u_2)_+\|_{L^\infty(B_{1/2}(x_0))}\leq C \int _{B_1(x_0)}(u_2)_+(x)dx \leq C \left(\int _{B_1(x_0)}(u_2)_+^\gamma(x)dx\right)^{\frac{1}{\gamma}}\leq C .$$
 Recalling that $u=u_1+u_2$ and the arbitrariness of $x_0$, we derive $u_+(x)\leq C$. Hence, $\sup\limits_{\mathbb{R}^n}u \leq C$.
$\hfill\square$

It's well known that if $f\in L^q({\mathbb{R}^n})$ and $g\in L^{q'}({\mathbb{R}^n})$ with $\frac{1}{q}+\frac{1}{q'}=1$, $q, q'>1$. Then the convolution $f\ast g\in C({\mathbb{R}^n})$ and  tends to 0 as $|x|\rightarrow \infty$. Combining  this, the result $\sup\limits_{\mathbb{R}^n}u_+\leq C$ and the integrability assumption of $u_+^\gamma$, we can obtain  the following limits.

\begin{lem}\label{L1}
Suppose $\gamma\in(1,\frac{n+2p}{n-2p})$, $n>2p$, $p\geq2$ and $u$ is a solution of equation (\ref{01}) satisfying $u(x)=o(|x|^2)$ at infinity. If $\zeta(x)$ is given by (\ref{ppffaa1}), then
$$
\lim\limits_{|x|\rightarrow \infty} (-\Delta)^i \zeta(x)=0, ~~i=0, 1,..., p-1.
$$
\end{lem}

\begin{lem}\label{a4}
Let $\gamma\in(1,\frac{n+2p}{n-2p})$, $n>2p$ and $p\geq 2$. Suppose $u$ is a nonconstant solution of (\ref{01}) satisfying satisfying $u(x)=o(|x|^2)$, then there exists $c_0<0$ such that
\begin{equation}\label{012}
u(x)=C_{n,p}\int_{\mathbb{R}^n}\frac{u_+^\gamma(y)}{|x-y|^{n-2p}}dy+c_0, ~~x\in{\mathbb{R}^n}.
\end{equation}
Moreover, the support of $u_+$ is compact.
\end{lem}

{\bf Proof } Lemma \ref{0723abc}, Lemma \ref{a2} and the fact that $\zeta\leq0$ give that
$$\Delta(u+\zeta)=0,~~~ u+\zeta\leq C ~~\mbox{ in } {\mathbb{R}^n}.$$
 Thus the Liouville theorem yields that there exists a constant $c_0$ such that
$$u(x)=-\zeta(x)+c_0.$$
We claim that $c_0< 0$. If $c_0\geq0$,  we have $u(x)\geq 0$ in ${\mathbb{R}^n}$, which is impossible from  the nonexistence  of positive solutions in subcritical exponent case of \cite{WJX}. Hence (\ref{012}) holds.

From (\ref{012}) and Lemma \ref{L1}, we obtain that $u$ tends to the negative constant $c_0$ as $|x|\rightarrow \infty$. As a consequence $u_+$ has a compact support.
$\hfill\square$

By applying the asymptotic expansions at infinity of the functions $g_i(z)=\frac{1}{|z|^{n-2p+2i}} (i=0, 1, \ldots, p-1 )$
\begin{align*}
\begin{split}
\frac{1}{|x-y|^{n-2p+2i}}&=\frac{1}{|x|^{n-2p+2i}}+(n-2p+2i)\sum_{j=1}^n\frac{x_jy_j}{|x|^{n-2p+2i+2}}
\\
&+\frac{(n-2p+2i)(n-2p+2i+2)}{2}\sum_{\ell,j=1}^n\frac{x_\ell x_j y_\ell y_j}{|x|^{n-2p+2i+4}}-\frac{(n-2p+2i)}{2}\sum_{j=1}^n\frac{y^2_j}{|x|^{n-2p+2i+2}}
\\
&
+\sum_{\ell,j,m=1}^ny_\ell y_j y_m\int_0^1\frac{(1-t)^2}{2}\partial_{\ell j m}g_i(ty)dt,
\end{split}
\end{align*}
and the fact that $u_+$ has a compact support, almost similar argument of Proposition 3.4 in \cite{RHA} gives the following asymptotic behavior of $\zeta$ and its derivatives.
\begin{lem}\label{L1999912}
Suppose $\gamma\in(1,\frac{n+2p}{n-2p})$, $n>2p$, $p\geq2$ and $u$ is a nonconstant solution of (\ref{01}) satisfying $u(x)=o(|x|^2)$. If $\zeta(x)$ is given by (\ref{ppffaa1}), then
$$
\lim\limits_{|x|\rightarrow \infty} (-\Delta)^i \zeta(x)|x|^{n-2p+2i}=a_i, ~~i=0, 1,..., p-1,
$$
where $a_i$ is given by
$$a_0=-C_{n,p}\int_{\mathbb{R}^n}u_+^\gamma(y)dy,~~a_{i+1}=(2p-2i-2)(n+2i-2p)a_i.$$
\end{lem}

Note that $a_i\leq0(i=1, \ldots, p-1)$.
From Lemma \ref{0723abc} and Lemma \ref{L1999912}, we have
\begin{align}\label{y1}
\lim\limits_{|x|\rightarrow \infty} (-\Delta)^i u(x)|x|^{n-2p+2i}=-a_i, ~~i=1,2,..., p-1.
\end{align}
Particularly, the result of (\ref{y1}) for the case $i=1$ yields that $\Delta u(x)\rightarrow 0$ as $|x|\rightarrow \infty$.

Now we will carry out the method of moving planes for integral equations (see \cite{11}) to complete the proof of Theorem \ref{T1}.

{\bf Proof of Theorem \ref{T1}} From Lemmas \ref{L1}-\ref{a4}, we obtain that
\begin{equation}\label{x2}
\zeta(x)=-C_{n,p}\int_{\mathbb{R}^n}\frac{(c_0-\zeta)_+^\gamma(y)}{|x-y|^{n-2p}}dy, ~~x\in{\mathbb{R}^n}
\end{equation}
and $\lim\limits_{|x|\rightarrow \infty}\zeta(x)=0$. To complete the proof of Theorem \ref{T1}, it is enough to show that $\zeta$ is symmetric about some point $x_0\in{\mathbb{R}^n}$ and $\frac{\partial \zeta}{\partial r}>0$ where $r=|x-x_0|$.

For $x=(x_1,x_2,...,x_n)\in {\mathbb{R}}^n$ and $\lambda\in{\mathbb{R}}$, we define ${\mathrm{T}}_{\lambda}=\{x\in{\mathbb{R}}^n|x_1=\lambda\}$, $\Sigma_\lambda=\{x\in{\mathbb{R}}^n|x_1<\lambda\}$, $x^\lambda=(2\lambda-x_1,x_2,...,x_n)$ and $\zeta_\lambda(x)=\zeta(2\lambda-x_1,x_2,...,x_n)=\zeta(x^\lambda)$.
Set $w_\lambda(x)=\zeta(x)-\zeta_\lambda(x)$. It's obvious that
\begin{equation*}
\zeta_\lambda(x)=-C_{n,p}\int_{\mathbb{R}^n}\frac{(c_0-\zeta_\lambda)_+^\gamma(y)}{|x-y|^{n-2p}}dy, ~~x\in{\mathbb{R}^n}.
\end{equation*}
From this and (\ref{x2}), we have
\begin{align}\label{x5}
\zeta_\lambda(x)-\zeta(x)=C_{n,p}\int_{\Sigma_\lambda}\left(\frac{1}{|x-y|^{n-2p}}-\frac{1}{|x-y^\lambda|^{n-2p}}\right)\left((c_0-\zeta)_+^\gamma(y)-(c_0-\zeta_\lambda)_+^\gamma(y)\right)dy.
\end{align}
\par $\mathbf{Step~ 1}$: We claim that for $\lambda$ sufficiently negative,
\begin{align}\label{x3}
w_\lambda(x)> 0, ~~~~~x\in \Sigma_\lambda.
\end{align}
Due to $\lim\limits_{|x|\rightarrow \infty}\zeta(x)=0$ and $c_0<0$, we have for $\lambda$ sufficiently negative
$$ (c_0-\zeta)_+^\gamma-(c_0 -\zeta_\lambda)_+^\gamma=-(c_0 -\zeta_\lambda)_+^\gamma\leq 0,~~~~x\in \Sigma_\lambda.$$
From this, (\ref{x5}) and the fact that $u$ is a nonconstant solution of (\ref{01}), we have $\zeta_\lambda(x)-\zeta(x)< 0$ for any $x\in \Sigma_\lambda$. Thus (\ref{x3}) holds.

\par $\mathbf{Step~ 2}$: Step 1 provides a starting point, from which we can now move the plane $T_\lambda$ to the right as long as (\ref{x3}) holds to its limiting position. Let
\begin{equation*}
\lambda_0=\sup\{\lambda|w_\mu(x)>0, \forall x \in \Sigma_\mu, \mu\leq\lambda\}.
\end{equation*}
It's clearly that $\lambda_0<+\infty$ and
$$ w_{\lambda_0}(x)\geq 0, ~~~x\in\Sigma_{\lambda_0}.$$
We will show that $ w_{\lambda_0}(x)\equiv0$ for $x\in\Sigma_{\lambda_0}$.

Otherwise if $w_{\lambda_0}\geq 0$ and $w_{\lambda_0}\not\equiv0$, we must have
\begin{align}\label{x7}
w_{\lambda_0}(x)>0,~~~x\in\Sigma_{\lambda_0},
\end{align}
where (\ref{x7}) follows from (\ref{x5}). From the facts  $\lim\limits_{|x|\rightarrow \infty}\zeta(x)=0$ and $c_0<0$, we derive that there exists $R$ large enough such that
\begin{equation}\label{p1}
(c_0-\zeta)_+\equiv0, ~~x\in\mathbb{R}^n\setminus B_R(0).
\end{equation}
Fixing this $R$, we have there exists constant $\delta>0$ and $c>0$ such that
\begin{align}\label{x8}
w_{\lambda_0}(x)\geq c,   ~~~x\in\overline{\Sigma_{\lambda_0 -\delta}\cap B_R(0)}.
\end{align}
Therefore by the continuity of $w_\lambda$ in $\lambda$ there exists $\varepsilon>0$ and $\varepsilon<\delta$ such that for all $\lambda\in[\lambda_0,\lambda_0+\varepsilon)$, we have
\begin{equation*}
w_\lambda(x)\geq0, ~~~x\in\overline{\Sigma_{\lambda_0 -\delta}\cap B_R(0)}.
\end{equation*}
We will show that for sufficiently small $0<\varepsilon<\delta$ and any $\lambda\in[\lambda_0,\lambda_0+\varepsilon)$
\begin{align}\label{x11}
w_\lambda(x)\geq0, ~~~x\in\Sigma_\lambda,
\end{align}
which contradicts with the definition of $\lambda_0$. Therefore we must have $w_{\lambda_0}\equiv 0$. Define
$$\Sigma_\lambda^-=\{x\in\Sigma_\lambda|w_\lambda(x)<0\}.$$
Next we claim that $\Sigma_\lambda^-$ must be measure zero.

For $y\in\Sigma_\lambda^-$, we can obtain that
\begin{equation}\label{4.1}
(c_0 -\zeta)_+^\gamma(y)-(c_0-\zeta_\lambda)_+^\gamma(y)\leq \gamma (c_0 -\zeta)_+^{\gamma-1}(y)|w_\lambda(y)|.
\end{equation}
Thus for $x\in \Sigma_\lambda$,
\begin{align}\label{x22}
\begin{split}
\zeta_\lambda(x)-\zeta(x)&\leq C_{n,p}\int_{\Sigma_\lambda^-}\left(\frac{1}{|x-y|^{n-2p}}-\frac{1}{|x-y^\lambda|^{n-2p}}\right)\left((c_0- \zeta)_+^\gamma(y)-(c_0-\zeta_\lambda)_+^\gamma(y)\right)dy \\
&\leq C_{n,p}\int_{\Sigma_\lambda^-}\left(\frac{1}{|x-y|^{n-2p}}-\frac{1}{|x-y^\lambda|^{n-2p}}\right)\gamma (c_0-\zeta(y))_+^{\gamma-1}|w_\lambda(y)|dy.
\end{split}
\end{align}
Applying Hardy-Littlewood-Sobolev inequality \cite{Dai} and H\"{o}lder inequality to (\ref{x22}) we obtain that
\begin{align}\label{x21}
\begin{split}
\|\zeta_\lambda(x)-\zeta(x)\|_{L^{\frac{2n}{n-2p}}(\Sigma_\lambda^-)}&\leq C\left(\int_{\Sigma_\lambda^-}\left((c_0-\zeta)_+^{\gamma-1}(y)|w_{\lambda}(y)|\right)^{\frac{2n}{n+2p}}dy\right)^{\frac{n+2p}{2n}}\\
&\leq C \left(\int_{\Sigma_\lambda^-}\left((c_0-\zeta)_+^{\gamma-1}(y)\right)^{\frac{n}{2p}}dy\right)^{\frac{2p}{n}}\left(\int_{\Sigma_\lambda^-}|w_\lambda(y)|^{\frac{2n}{n-2p}}dy\right)^{\frac{n-2p}{2n}}.
\end{split}
\end{align}
Recall that $\Sigma_\lambda^-\subset ((\Sigma_\lambda\setminus \Sigma_{\lambda_0 -\delta})\cap B_R)\cup (\Sigma_\lambda \setminus  B_R)$ and $-\zeta$ is bounded above, we can choose $\delta$ sufficiently small such that
$$ C \left(\int_{{\Sigma_\lambda^-}\cap B_R}\left((c_0-\zeta)_+^{\gamma-1}(y)\right)^{\frac{n}{2p}}dy\right)^{\frac{2p}{n}}\leq \frac{1}{2}. $$
From this and (\ref{p1}), we have
$$ C \left(\int_{\Sigma_\lambda^-}\left((c_0-\zeta)_+^{\gamma-1}(y)\right)^{\frac{n}{2p}}dy\right)^{\frac{2p}{n}}\leq \frac{1}{2}. $$
Now (\ref{x21}) implies that $\|w_\lambda\|_{L^{\frac{2n}{n-2p}}(\Sigma_\lambda^-)}=0$ and therefore $\Sigma_\lambda^-$ must be measure zero.

This verifies (\ref{x11}). Thus we must have $w_{\lambda_0}\equiv0$.

\par $\mathbf{Step~ 3}$: We show that $\frac{\partial \zeta}{\partial x_1}<0$ for $x\in\Sigma_{\lambda_0}$.

In fact, from the definition of $\lambda_0$ we have for any $\lambda<\lambda_0$,
\begin{equation}\label{z1}
w_\lambda(x)> 0, ~~~~~x\in \Sigma_\lambda.
\end{equation}
Simple calculation gives that for any $x\in T_\lambda$ with $\lambda<\lambda_0$,
\begin{equation*}
\begin{split}
\zeta_{x_1}(x)=&C_{n,p}(n-2p)\int_{\mathbb{R}^n}\frac{(c_0-\zeta)_+^\gamma(y)(x_1-y_1)}{|x-y|^{n-2p+1}}dy\\
=&C_{n,p}(n-2p)\int_{\Sigma_\lambda}\frac{(c_0-\zeta)_+^\gamma(y)(x_1-y_1)}{|x-y|^{n-2p+1}}dy+C_{n,p}(n-2p)\int_{\Sigma_\lambda}\frac{(c_0-\zeta_\lambda)_+^\gamma(y)(x_1-y_1^\lambda)}{|x-y^\lambda|^{n-2p+1}}dy\\
=&C_{n,p}(n-2p)\int_{\Sigma_\lambda}\frac{\left((c_0-\zeta)_+^\gamma(y)-(c_0-\zeta_\lambda)_+^\gamma(y)\right)(x_1-y_1)}{|x-y|^{n-2p+1}}dy\\
<&0,
\end{split}
\end{equation*}
where the last inequality follows from (\ref{z1}). Thus the claim holds.

Since the problem is invariant with respect to rotation, we can take any direction as the $x_1$ direction. Hence we have that $\zeta$ is radially symmetric about some $x_0\in\mathbb{R}^n$ and $\frac{\partial \zeta}{\partial r}>0$ where $r=|x-x_0|$.
$\hfill\square$

Actually, we may also prove Theorem \ref{T1}  by applying  moving plane method to the function $(-\Delta)^{p-1}u$, after asymptotic behaviors at infinity of this function and its first-order derivatives are established. Similar argument can be seen  in \cite{RHA}.

\subsection{Decay assumption $\Delta u(x)\rightarrow 0$ as $|x|\rightarrow \infty$}$\\$

Assume that $\Delta u(x)\rightarrow 0$ as $|x|\rightarrow \infty$, then $\Delta u\in L^\infty(\mathbb{R}^n)$.
From the proof of Lemmas \ref{a2} and \ref{L1}, we know that the results of these two lemmas are still valid replacing the growth assumption $u(x)=\mbox{o}(|x|^2)$ at infinity by the decay assumption  $\Delta u(x)\rightarrow 0$ as $|x|\rightarrow \infty$.

If we can prove that the results of Lemma \ref{a4}  still hold true under the  assumption that $\Delta u(x)\rightarrow 0$ as $|x|\rightarrow \infty$, instead of the assumption $u(x)=\mbox{o}(|x|^2)$ at infinity, then we obtain that $u(x)\rightarrow c_0<0$ as $|x|\rightarrow \infty$, which gives that $u(x)=\mbox{o}(|x|^2)$ at infinity. Then Theorem \ref{T111} follows from  Theorem \ref{T1} immediately.

To obtain the results of Lemma \ref{a4} (with the decay assumption that $\Delta u(x)\rightarrow 0$ as $|x|\rightarrow \infty$), from the proof of Lemma \ref{a4}, it is enough to prove Lemma \ref{0723}. To this end, we need to establish the  super polyharmonic properties, namely Lemma \ref{a1220727}.

\begin{lem}\label{a1220727}
Suppose that $\gamma\in(1,\frac{n+2p}{n-2p})$, $n>2p$, $p\geq 2$ and $u$ is a solution of equation (\ref{01}) satisfying $\Delta u(x)\rightarrow 0$ as $|x|\rightarrow \infty$, then we have $$(-\Delta)^i u\geq 0, ~~~i=1,2,...,p-1.$$
\end{lem}
{\bf Proof } The proof of this lemma is almost similar to Lemma \ref{a1}. For the convenience of readers, we will briefly introduce the proof here. Denote $w_i=(-\Delta)^iu$. We  claim that $w_i \geq 0$, $i=2,...,p$.

Obviously this claim is  true when $p=2$. For the cases $p>2$, we have $w_p\geq 0$ due to $(-\Delta)^p u(x)=u_+^\gamma(x)\geq 0$ in $\mathbb{R}^n$. We next prove that $w_{p-1}\geq 0$ by contradiction. Suppose there exists $x_0 \in {\mathbb{R}^n}$ such that $w_{p-1}(x_0)< 0$. Without loss of generality, we assume that $x_0=0$. Arguments similar to that of Lemma \ref{a1} show that $(-1)^i\bar{w}_{p-i}(r)\geq c_i r^{2(i-1)}$ for $ r\geq r_{i-1}, i=1,...,p-1$. Particularly, we have
$$ (-1)^{p-1}\bar{w}_1(r)\geq c_{p-1} r^{2(p-2)}~~~~\mbox{ for } r\geq r_{p-2}.   $$
Hence $|\bar{w}_1(r)|\geq cr^2$ for $r\geq r_{p-2}$ due to $p>2$, which is a contradiction. Thus $w_{p-1}\geq 0$. Similar induction process  as that of Lemma \ref{a1} shows that $w_{p-i} \geq 0$, $i=2,...,p-2$.
Note that $i$ cannot arrive at $p-1$, since our assumption is $\Delta u(x)\rightarrow 0$ as $|x|\rightarrow \infty$ rather that the assumption $u(x)=\mbox{o}(|x|^2)$ at infinity as that of Lemma \ref{a1}.

Especially we have $(-\Delta)^2u (x)\geq 0$ in $\mathbb{R}^n$. From this and the assumption $\Delta u(x)\rightarrow 0$ as $|x|\rightarrow \infty$, applying the maximum principle we have that $-\Delta u(x)\geq 0$ in $\mathbb{R}^n$. Hence we complete the proof of this lemma.
$\hfill\square$

\begin{lem}\label{0723}
Suppose that $\gamma\in(1,\frac{n+2p}{n-2p})$, $n>2p$, $p\geq 2$ and $u$ is a solution of equation (\ref{01}) satisfying $\Delta u(x)\rightarrow 0$ as $|x|\rightarrow \infty$. Let $\zeta(x)$ be given by (\ref{ppffaa1}). Then $(-\Delta)^i(u+\zeta)=0$, $i=1,2,...,p-1$.
\end{lem}

{\bf Proof} For the case of $i=p-1$, set $w(x)=(-\Delta)^{p-1}(u+\zeta)(x)$. Then $\Delta w(x)=0$. From Lemma \ref{a1220727} and the results of Lemma \ref{L1}  (with the decay assumption that $\Delta u(x)\rightarrow 0$ as $|x|\rightarrow \infty$),  we have $\lim\limits_{|x|\rightarrow \infty}w(x)\geq0$. Then Liouville theorem implies that there exists $C\geq0$ such that $w(x)\equiv C$.
\par Next, we claim that $C=0$. If not, we have $(-\Delta)^{p-1}u(x)\geq \frac{C}{2}$ for $|x|$ large. Let $\bar{u}$ be the average of $u$ which is defined as before. Then as in Lemma \ref{a1220727}, we have
$$ (-1)^{p-2}(-\Delta\bar{u}(r))\geq C_0 r^{2(p-2)},$$
for $r$ large enough. Clearly it contradicts with the assumption that $\Delta u(x)\rightarrow 0$ as $|x|\rightarrow \infty$, since $p\geq2$. So we have $w(x)\equiv0$.

Similar arguments as the above, we obtain $(-\Delta)^i(u+\zeta)=0$, $i=1,...,p-2$.
$\hfill\square$

\section{Radial symmetry of solutions to fractional  equations with integral constraint}
\par For the convenience of readers, we first recall regularity results and Harnack inequality of linear fractional Laplace equations that will be used in the forthcoming sections.
\subsection{Regularity results and nonlocal Harnack inequality}

\par Ros-Oton and Serra in \cite{RSL} established the Schauder interior estimates of Poisson equation involving the fractional Laplacian, which is a counterpart of  classical Poisson equation. After scaling, one can obtain the following two propositions. Denote
$$[u]_{C^{k,\alpha}(B_R)}=\sum\limits_{|\eta|=k}\sup\limits_{{x,y\in B_R} \atop {x\neq y}}\frac{|D^\eta u(x)-D^\eta u(y)|}{|x-y|^\alpha}.$$

\begin{pro}\label{FL1}
  Assume that $f\in L^\infty (B_{2R})$ and $u\in L^\infty (B_{2R})\cap \mathcal{L}_{2s}({\mathbb{R}}^n)$ is a solution of $(-\Delta)^su=f(x)$ in $B_{2R}$, then $u\in C^\alpha(B_{R})$ for any $\alpha\in (0,\min\{2s,1\})$ and
\begin{equation*}
\begin{split}
\|u\|_{L^\infty(B_R)}+R^{\alpha}[u]_{C^{0,\alpha}(B_R)}\leq  C [R^{2s}\|f\|_{L^\infty(B_{2R})}+\frac{1}{R^n}\int_{{\mathbb{R}}^n}\frac{|u(y)|}{(1+|\frac{y}{R}|)^{n+2s}}dy+\|u\|_{L^\infty(B_{2R})}].\hspace{4.6cm}
\end{split}
\end{equation*}

\end{pro}

\begin{pro}\label{FL2}
 Suppose that $\alpha+2s= k+\beta$, where $0<\alpha, \beta<1$ and $k\in\mathbb{N}
 +$. If $f\in C^\alpha (B_{2R})$ and $u\in C^\alpha (B_{2R})\cap \mathcal{L}_{2s}({\mathbb{R}^n})$ is a solution of
$ (-\Delta)^su(x)=f(x)$ in $B_{2R}$, then
\begin{equation*}
\begin{split}
\sum_{|\eta|=1}^{k}R^{|\eta|\ }|D^\eta u\|_{L^\infty(B_R)}+R^{2s+\alpha}[u]_{C^{k,\beta}(B_R)}\hspace{9.4cm}\\
\leq  C [R^{2s}\|f\|_{L^\infty(B_{2R})}+R^{2s+\alpha}[f]_{C^{0,\alpha}(B_{2R})}+\frac{1}{R^n}
\int_{{\mathbb{R}}^n}\frac{|u(y)|}{(1+|\frac{y}{R}|)^{n+2s}}dy+\|u\|_{L^\infty(B_{2R})}+R^{\alpha}[u]_{C^{0,\alpha}(B_{2R})}].
\end{split}
\end{equation*}
\end{pro}

\par The following is the nonlocal Harnack inequality, given in \cite{K.M}.
\begin{pro}\label{FL3}
 (Nonlocal Harnack inequality) Suppose that $u\in C_{\mbox{loc}}^{2s+\alpha} \cap \mathcal{L}_{2s}({\mathbb{R}}^n)$ satisfies
\begin{eqnarray}
\nonumber
\begin{cases}
(-\Delta)^su=0  & \mbox{ in }{B_R},~~~~~~\\
u=g~~~~~~~&\mbox{ in } ~{\mathbb{R}^n}\backslash {B_R}
\end{cases}
\end{eqnarray}
with $u\geq0$ in $B_R$, then  the following estimate holds
$$ \sup_{B_{\frac{R}{2}}}u\leq c \inf_{B_{\frac{R}{2}}}u+c\int_{{\mathbb{R}^n}\backslash {B_R}}\frac{g_{-}(z)}{|z|^{n+2s}}dz,$$
where the positive constant $c$ depends only on $n,s,R$.

\end{pro}

\subsection{ Various maximum principles }
\par We introduce several maximum principles: Propositions \ref{FL4}-\ref{FL7}, established in \cite{chen}, which will play a key role in  using of the method of moving planes.
\begin{pro}\label{FL4}
 (Maximum principle). Let $\Omega$ be a bounded domain in ${\mathbb{R}^n}$. Assume that $u \in \mathcal{L}_{2s} \cap C_{\mbox{loc}}^{1,1}(\Omega)$ and is lower semi-continuous on $\bar{\Omega}$. If
\begin{equation*}
\left\{\begin{array}{rl}
(-\Delta)^su(x)\geq 0  ~~  &\mbox{in}~ \Omega,~~~~~\\
u(x)\geq0~~~~~~~~&\mbox{in} ~ {\mathbb{R}^n}\backslash \Omega.
\end{array}\right.
\end{equation*}
then $u(x)\geq0~~\mbox{in} ~\Omega.$
\end{pro}

\begin{pro}\label{FL5}
(Maximum Principle for Anti-symmetric Functions) Let $\mathrm{T}_\lambda$ be a hyperplane in ${\mathbb{R}^n}$. Without loss of generality, we may assume that $\mathrm{T}_\lambda=\{x\in {\mathbb{R}^n}|x_1=\lambda, \mbox{for some } \lambda \in {\mathbb{R}}\}$. Denote $ {\tilde{x}}=(2\lambda-x_1,x_2,...,x_n)$ and $ {\Sigma_\lambda}=\{x\in {\mathbb{R}^n}|x_1<\lambda\}$.
Let $\Omega$ be a bounded domain in $\Sigma_\lambda$. Assume that $u\in \mathcal{L}_{2s}  \cap C_{\mbox{loc}}^{1,1}(\Omega)$ and is lower semi-continuous on $\bar{\Omega}$. If
\begin{equation*}
\left\{\begin{array}{rl}
(-\Delta)^su(x)\geq 0  ~~  &\mbox{in}~ \Omega,~~~~~~\\
u(x)\geq0~~~~~~~~~~&\mbox{in} ~ \Sigma_\lambda\backslash \Omega,~\\
u(\tilde{x})=-u(x)~~~~&\mbox{in}~\Sigma_\lambda,~~~~
\end{array}\right.
\end{equation*}
then $u(x)\geq0$ in $\Omega$.

\end{pro}

\begin{pro}\label{FL6}
(Narrow Region Principle) Suppose that the definitions of $\mathrm{T}_\lambda$, ${\tilde{x}}$ and ${\Sigma_\lambda}$ are the same as in proposition \ref{FL5}. Let $\Omega$ be a bounded domain in $\Sigma_\lambda$, such that it is contained in $\{x|\lambda-l<x_1<\lambda\}$ with small $l$. Further assume that $u\in \mathcal{L}_{2s}  \cap C_{\mbox{loc}}^{1,1}(\Omega)$ and is lower semi-continuous on $\bar{\Omega}$. If $c(x)$ is bounded from below in $\Omega$ and
\begin{equation*}
\left\{\begin{array}{rl}
(-\Delta)^su(x)+c(x)u(x)\geq 0  ~~  &\mbox{in}~ \Omega,~~~~~\\
u(x)\geq0~~~~~~~~~~~~~~~~~~~~~~&\mbox{in} ~ \Sigma_\lambda\backslash \Omega,~\\
u(\tilde{x})=-u(x)~~~~~~~~~~~~~~~~&\mbox{in}~\Sigma_\lambda,~~~~
\end{array}\right.
\end{equation*}
then for sufficiently small $l$, we have $u(x)\geq0$ in $\Omega$.
\end{pro}

 Moreover,  for any  $u$ of the above three propositions, if $u=0$ at some point in $\Omega$, then $ u\equiv0 $ in $ {\mathbb{R}^n}$. This and those conclusions in the above three propositions all  hold for unbounded region $\Omega$ if we further assume that $ \liminf\limits_{|x|\rightarrow \infty} u(x)\geq 0$.

\begin{pro}\label{FL7}
(Decay at infinity) Let $ {\Sigma_\lambda}=\{x\in {\mathbb{R}^n}|x_1<\lambda \mbox{ for some } \lambda \in {\mathbb{R}}\}$ and $\Omega$ be an unbounded region in ${\Sigma_\lambda}$. Assume $u \in \mathcal{L}_{2s}\cap C_{\mbox{loc}}^{1,1}(\Omega)$ is a solution of
\begin{eqnarray}\label{3.1}
\begin{cases}
(-\Delta)^su(x)+c(x)u(x)\geq 0  &  \mbox{in}~ \Omega,\\
u(x)\geq0~&\mbox{in}~ \Sigma_\lambda\backslash \Omega,\\
u(\tilde{x})=-u(x) &\mbox{in}~\Sigma_\lambda
\end{cases}
\end{eqnarray}
with
\begin{equation}\label{kkk}
\liminf\limits_{|x|\rightarrow \infty}|x|^{2s} c(x)\geq 0,
\end{equation}
then there exists a constant $R_0>0$ ( depending on $c(x)$, but is independent of $u$) such that if
\begin{equation}\label{oo}
u(x^0)=\min_{\Omega}u(x)<0,
\end{equation}
then
\begin{equation}\label{3.4}
|x^0|\leq R_0.
\end{equation}

\end{pro}
\par From the proof of Proposition \ref{FL7}, we know that the  same result holds true with the assumption $\liminf\limits_{|x\rightarrow \infty}|x|^{2s}c(x)\geq 0$ replaced by $c(x)>-\tilde{c}|x|^{-2s}$ for $|x|$ large enough, where $\tilde{c}=\frac{2w_n C_{s,n}}{4^{n+2s}}$.

\begin{pro*}\label{lemma* 7}
(Decay at infinity) Let $ {\Sigma_\lambda}=\{x\in {\mathbb{R}^n}|x_1<\lambda \mbox{ for some } \lambda \in {\mathbb{R}}\}$ and $\Omega$ be an unbounded region in ${\Sigma_\lambda}$. Assume $u \in \mathcal{L}_{2s}\cap C_{\mbox{loc}}^{1,1}(\Omega)$ is a solution of
(\ref{3.1})
with
\begin{equation}\label{3.2}
c(x)>-\tilde{c}|x|^{-2s}~~\mbox{ for } |x| >R_0,
\end{equation}
where $\tilde{c}=\frac{2w_n C_{s,n}}{4^{n+2s}}$ and $R_0>0$. If
\begin{equation}\label{3.3}
u(x^0)=\min_{\Omega}u(x)<0,
\end{equation}
then
\begin{equation}\label{3.4}
|x^0|\leq \max\{R_0, |\lambda|\}.
\end{equation}

\end{pro*}

\begin{rem}\label{FR2}
\rm Note that results in Propositions \ref{FL4}-$3.7^{(')}$ are also true if $u \in \mathcal{L}_{2s}\cap C_{\mbox{loc}}^{2s+\alpha}(\Omega)$,  even for the case $2s+\alpha<2$.  Since the only role of the regularity $u \in \mathcal{L}_{2s}\cap C_{\mbox{loc}}^{1,1}(\Omega)$ played in \cite{chen} is to make the operator $(-\Delta)^su$ well defined. Recall that $(-\Delta)^s$ is well defined for functions belonging to $\mathcal{L}_{2s}\cap C_{\mbox{loc}}^{2s+\alpha}(\Omega)$.
\end{rem}

\subsection{ Proof of Theorem \ref{F1} }

\par Next, we will prove Theorem \ref{F1}. To this end, we need to establish the following several lemmas.

\begin{lem}\label{FP1}
 Given $R>0$ and $M>0$. Let $v=v(x)\in \mathcal{L}_{2s}\cap C_{\mbox{loc}}^{2s+\alpha}(\mathbb{R}^n)$ be a solution of
\begin{equation}\label{x40}
\left\{\begin{array}{rl}
(-\Delta)^sv=v^\gamma_+  ~~ &\mbox{ in }{\mathbb{R}}^n,\hspace{2.1cm}~\\
v\leq M~~~~~~~~~ &\mbox{ in }B_R,~\hspace{2cm}\\
v(x_0)=1 ~~~~~~~ &\mbox{for some} ~~ x_0 \in B_{\frac{R}{2}}.
\end{array}\right.
\end{equation}
Then there exists $C_0=C_0(n, s, \gamma, R, M, \|v\|_{\mathcal{L}_{2s}})>0$ such that
\begin{equation*}
v\geq-C_0 ~~~~~  \mbox{in}  ~B_{\frac{R}{4}}.
\end{equation*}

\end{lem}
{\bf Proof }
Denote $v_1, v_2$ as the solutions of the following problems respectively
\begin{eqnarray} \nonumber
\begin{cases}
(-\Delta)^sv_1=v_+^\gamma  &\mbox{ in}~ B_R,\\
v_1=0 &\mbox{in}~{\mathbb{R}}^n\backslash B_R,
\end{cases}
\indent
\begin{cases}
(-\Delta)^sv_2=0 & \mbox{in} ~B_R,\\
v_2=v &\mbox{in}~{\mathbb{R}}^n\backslash B_R.
\end{cases}
\end{eqnarray}
Then their sum $w=v_1+v_2$ satisfies
\begin{eqnarray} \nonumber
\begin{cases}
(-\Delta)^sw=v^\gamma_+  & \mbox{ in }{B_R},\\
w=v  &\mbox{ in } {\mathbb{R}}^n \setminus {B_R}.
\end{cases}
\end{eqnarray}
We have $v\equiv w$ in ${\mathbb{R}}^n$, since Proposition \ref{FL4} tells us that trivial solution is the unique solution of the problem
$(-\Delta)^s\phi=0   \mbox{ in }{B_R},~
\phi=0 ~\mbox{ in } {\mathbb{R}}^n \setminus {B_R}$.
So $v=v_1+v_2$. Kulczycki in \cite{Ku} showed  that the Green's function on the ball $B_R(0)$ corresponding to $(-\Delta)^s$ with homogeneous Dirichlet boundary condition is
\begin{equation*}
G_R(x,y)={\frac{A_{n,s}}{{\tau_R}^{\frac{n-2s}{2}}}}\left[1-\frac{B_{n,s}}{(\tau_R+t_R)^{\frac{(n-2s)}{2}}} {\int_0^{\frac{\tau_R}{t_R}}} {\frac{(\tau_R-t_Rb)^{\frac{(n-2s)}{2}}}{b^s(1+b)}}db\right], ~\forall ~x,y\in B_R(0),
\end{equation*}
where $\tau_R=\frac{|x-y|^2}{R^2}$, $t_R=(1-\frac{|x|^2}{R^2})(1-\frac{|y|^2}{R^2})$, $A_{n,s}$ and $B_{n,s}$ are constants depending on $n$ and $s$. Then
\begin{equation*}
v_1(x)=\int_{B_R}G_R(x,y)v_+^\gamma(y)dy,
\end{equation*}
and so
\begin{equation}\label{4.4}
|v_1|\leq M^\gamma \sup_{x\in B_R}\|G_R(x,.)\|_{L^1(B_R)}\equiv C_1(R)  ~~ \mbox{in} ~ B_R.
\end{equation}
The reason of the final equality is the integrability of Green's function, we refer to \cite{Bucur} for details. Moreover, it follows from the maximum principle that $v_1\geq 0$, and so $v_2=v-v_1\leq v$ in $B_R$. Then $v_2\leq M$ in $B_R$. We have
\begin{equation*}
\max_{{\bar{B}_{\frac{R}{2}}}}v_2\geq\max_{{\bar{B}_{\frac{R}{2}}}}v-\max_{{\bar{B}_{\frac{R}{2}}}}v_1\geq v(x_0)-\max_{{\bar{B}_{\frac{R}{2}}}}v_1\geq 1-C_1.
\end{equation*}
Applying Proposition \ref{FL3} (Nonlocal Harnack inequality) to the nonnegative function $M-v_2$, we deduce that there exists $C_2=C_2(n, \gamma, R, M, s)>0$ such that
\begin{equation*}
\begin{split}
\max_{{{\bar{B}_{\frac{R}{4}}}}}(M-v_2)&\leq C_2\min_{{{\bar{B}_{\frac{R}{2}}}}} (M-v_2)+C_2\int_{{\mathbb{R}^n}\backslash {B_{{R}}}}\frac{(M-v_2)_{-}(z)}{|z|^{n+2s}}dz~~~~~\\
&\leq C_2\{M-(1-C_1)\}+C_2\int_{{\mathbb{R}^n}\backslash {B_{{R}}}}\frac{(M-v_2)_{-}(z)}{|z|^{n+2s}}dz\\
&\leq C_2\{M-(1-C_1)\}+C_2\|v\|_{\mathcal{L}_{2s}}.\hspace{2.8cm}
\end{split}
\end{equation*}
Therefore
\begin{equation}\label{4.5}
\min_{{{\bar{B}_{\frac{R}{4}}}}}v_2\geq M-C_2\{M-(1-C_1)\}-C_2\|v\|_{\mathcal{L}_{2s}}.
\end{equation}
The result of this lemma follows from (\ref{4.4}), (\ref{4.5}) and the fact that $ v=v_1+v_2$.
$\hfill\square$

\begin{lem}\label{FP2}
 Any solution $v=v(x)\in \mathcal{L}_{2s}\cap C_{\mbox{loc}}^{2s+\alpha}(\mathbb{R}^n)$ of (\ref{1}) satisfying $|v(x)|= \mbox{O}(|x|^\tau)$ at infinity with $\tau<2s$ is bounded above.
\end{lem}
{\bf Proof } Let $v=v(x)$ be a classical solution of (\ref{1}). Since $\int_{\mathbb{R}^n}v_+^{\frac{n(\gamma-1)}{2s}}dx<+\infty$,
then for any $k\in \mathbb{N}$, there exists a sequence $\{R_k\}$ tending to $+\infty$ such that
$$\int_{|x|\geq R_k}v_+^{\frac{n(\gamma-1)}{2s}}dx\leq \frac{1}{k}.  $$
If $v$ is not bounded from above, then there exists a sequence of points $\{x_k\} \subset \mathbb{R}^n$ satisfying $|x_k|>R_k+2$ and
\begin{equation*}
\left\{\begin{array}{rl}
\int_{B_1(x_k)}v_+^{\frac{n(\gamma-1)}{2s}}dx\leq\frac{1}{k},\hspace{2.8cm}\\
\max\limits_{\bar{B}_{{\frac{1}{4}}}(x_k)}v\geq k.\hspace{4.4cm}
\end{array}\right.
\end{equation*}
We define
$$  v_k(x)=v(x+x_k)   ~~~\mbox{ in } \mathbb{R}^n,$$
then $v_k$ satisfies
\begin{equation}\label{220803}
\left\{\begin{array}{rl}
(-\Delta)^sv_k=({v_k})^\gamma_+  ~~~ \mbox{ in }{\mathbb{R}^n},\hspace{2cm}\\
\int_{B_1(0)}({v_k})_+^{\frac{n(\gamma-1)}{2s}}dx\leq\frac{1}{k},\hspace{2.6cm}\\
\max\limits_{\bar {B}_{{\frac{1}{4}}}(0)}{v_k}\geq k.\hspace{4.2cm}
\end{array}\right.
\end{equation}
\par For each $k$, we define $h_k\in C^{2s+\alpha}(B_1)$ and $y_k \in B_{\frac{1}{2}}(0)$ by
\begin{equation*}\label{4.9}
h_k(y)=(\frac{1}{2}-r)^q v_k(y),~~ h_k(y_k)=\max_{\bar{B}_{\frac{1}{2}}}h_k(y),
\end{equation*}
where $q=\frac{2s}{\gamma-1}$ and $r=|y|$. It holds that
\begin{equation}\label{4.10}
\begin{split}
h_k(y_k)\geq \max_{\bar{B}_{\frac{1}{4}}} (\frac{1}{2}-r)^q v_k(y)
\geq (\frac{1}{4})^q \max_{\bar{B}_{\frac{1}{4}}} v_k(y)\geq (\frac{1}{4})^q k,
\end{split}
\end{equation}
for any $k$. Set
\begin{equation*}
\sigma_k=\frac{1}{2}-|y_k|, ~~d_k^q=h_k(y_k),~~ \mu_k=\sigma_k /d_k.
\end{equation*}
Then $v_k(y_k)=\frac{d_k^q}{\sigma_k^q}=\mu_k^{-q}$. For any $y\in B_{\frac{\sigma_k}{2}}(y_k)\subset B_{\frac{1}{2}}(0)$, we can obtain that
\begin{equation*}
\frac{1}{2}-|y|\geq \frac{1}{2}-(|y_k|+|y-y_k|)= (\frac{1}{2}-|y_k|)-|y-y_k|\geq \sigma_k-\frac{\sigma_k}{2}=\frac{\sigma_k}{2},
\end{equation*}
thus
\begin{equation}\label{4.11}
d_k^q=h_k(y_k)\geq (\frac{1}{2}-|y|)^qv_k(y)\geq (\frac{\sigma_k}{2})^q v_k(y).
\end{equation}
We introduce
\begin{equation*}
w_k(y)=\mu _k^qv_k(\mu _k y+y_k).
\end{equation*}
 From (\ref{220803}) and (\ref{4.11}), we can obtain that
\begin{equation}\label{4.12}
\left\{\begin{array}{rl}
(-\Delta)^s{w_k}=(w_k)^\gamma_+,~ w_k\leq 2^q ~~ \mbox{ in }{B_{\frac{d_k}{2}}},\hspace{1.4cm}\\
\int_{B_{\frac{d_k}{2}}}(w_k)_+^{\frac{n(\gamma-1)}{2s}}=\int_{B_{\frac{\sigma_k}{2}}(y_k)}(v_k)_+^{\frac{n(\gamma-1)}{2s}}\leq\frac{1}{k},\hspace{0.6cm}\\
w_k(0)=\mu_k^q v_k(y_k)=1.\hspace{3.8cm}
\end{array}\right.
\end{equation}

We claim that there exists a positive constant $C$ independent of $k$ such that
\begin{equation}\label{asa}
\int_{{\mathbb{R}}^n}\frac{|w_k(x)|}{1+|x|^{n+2s}}dx\leq C.
\end{equation}
 Due to $|v(x)|=\mbox{O}(|x|^\tau)$ at infinity, there exists $R>0$ and $C_0>0$ such that
$\frac{C_0}{2}|x|^\tau \leq |v(x)| \leq 2 C_0 |x|^\tau \mbox{ for any } |x|\geq R.$  Note that $d_k \rightarrow \infty$ by (\ref{4.10}), hence $\frac{d_k}{2}\geq R$ for $k$ large enough. Combined this and (\ref{4.12}), we find that for $k$ large enough
\begin{equation}\label{mm}
\begin{split}
\int_{{\mathbb{R}}^n}\frac{|w_k(x)|}{1+|x|^{n+2s}}dx&=\int_{B_R}\frac{|w_k(x)|}{1+|x|^{n+2s}}dx+\int_{B_R^c}\frac{|w_k(x)|}{1+|x|^{n+2s}}dx\hspace{4.8cm}\\
&\leq \int_{B_R}\frac{2^q}{1+|x|^{n+2s}}dx+\int_{B_R^c}\frac{|w_k(x)|}{1+|x|^{n+2s}}dx.\hspace{4.7cm}\\
\end{split}
\end{equation}
Obviously, the first integral is finite and its value is independent of $k$. The second integral need more rigorous  analysis.
Indeed, there exists $C_1>0$ independent of $k$ such that
\begin{equation}\label{nn}
\begin{split}
&\int_{B_R^c}\frac{|w_k(x)|}{1+|x|^{n+2s}}dx\hspace{12.6cm}~~\\
&=\int_{B_R^c}\frac{|\mu _k^qv(\mu _k x+y_k+x_k)|}{1+|x|^{n+2s}}dx\hspace{10.8cm}~~\\
&=\int_{{B_R^c}\cap \{|\mu _k x+y_k+x_k|\leq R\}}\frac{|\mu _k^qv(\mu _k x+y_k+x_k)|}{1+|x|^{n+2s}}dx
+\int_{{B_R^c}\cap \{|\mu _k x+y_k+x_k|\geq R\}}\frac{|\mu _k^qv(\mu _k x+y_k+x_k)|}{1+|x|^{n+2s}}dx~~~\\
&\leq C_1+\int_{{B_R^c}\cap \{|\mu _k x+y_k+x_k|\geq R\}}\frac{|\mu _k^qv(\mu _k x+y_k+x_k)|}{1+|x|^{n+2s}}dx,\hspace{7.2cm}\\
\end{split}
\end{equation}
where we used the fact that $\mu_k\rightarrow 0$ and $v\in C_{\mbox{loc}}^{2s+\alpha}(\mathbb{R}^n)$. Set
$$
A_{R,k}={B_R^c}\cap \{|\mu _k x+y_k+x_k|\geq R\}\cap \{|x|\leq |x_k|\},
$$
$$
 B_{R,k}={B_R^c}\cap \{|\mu _k x+y_k+x_k|\geq R\}\cap \{|x|\geq |x_k|\}.
$$
It follows from $\mu_k\rightarrow 0$ as $k\rightarrow \infty$ that $\mu_k<\frac{1}{2}$ for $k$ large enough. Combined this with the facts $|y_k|\leq \frac{1}{2}$ and $|x_k|\rightarrow \infty$, we conclude that for $x\in A_{R,k}$ and large $k$,
$$R\leq\frac{1}{2}|x_k| \leq|\mu _k x+y_k+x_k|\leq \frac{3}{2}|x_k|.$$
Therefore there exists $C_2>0$ ( independent of $k$ ) such that for large $k$,
\begin{equation}\label{bb}
\begin{split}
\int_{A_{R,k}}\frac{|\mu _k^qv(\mu _k x+y_k+x_k)|}{1+|x|^{n+2s}}dx&=\int_{A_{R,k}}\frac{|v(\mu _k x+y_k+x_k)|}{v(y_k+x_k)(1+|x|^{n+2s})}dx\hspace{2.6cm}\\
&\leq \int_{A_{R,k}}\frac{2C_0(\frac{3}{2}|x_k|)^\tau}{(\frac{C_0}{2})(\frac{1}{2}|x_k|)^\tau(1+|x|^{n+2s})}dx\hspace{2.3cm}\\
&\leq C_2.\hspace{7.4cm}
\end{split}
\end{equation}
Similarly, for any $x\in B_{R,k}$ and large $k$, we can obtain that  $|\mu_k x+y_k+x_k|\leq 2|x|$. Hence there exists $C_3>0$ ( independent of $k$ ) such that for large $k$
\begin{equation}\label{vv}
\begin{split}
\int_{B_{R,k}}\frac{|\mu _k^qv(\mu _k x+y_k+x_k)|}{1+|x|^{n+2s}}dx&\leq \int_{B_{R,k}}\frac{\mu _k^q\cdot2C_0\cdot(|\mu _k x+y_k+x_k|)^\tau}{1+|x|^{n+2s}}dx\hspace{2.2cm}\\
&\leq \int_{B_{R,k}}\frac{\mu _k^q\cdot2C_0\cdot 2^\tau\cdot|x|^\tau}{1+|x|^{n+2s}}dx\hspace{4cm}\\
&\leq C_3,\hspace{7.8cm}
\end{split}
\end{equation}
where we used the assumption that $\tau <2s$. The estimate (\ref{mm})-(\ref{vv}) implies that (\ref{asa}) holds.

For any given $R>4$, from (\ref{4.12}), (\ref{asa}) and Lemma \ref{FP1},  the uniformly boundedness of $w_k$  in $B_R$ follows. Combined this, (\ref{asa})
and  Proposition \ref{FL1}, we know that $\|w_k\|_{C^\alpha}$ is uniformly bounded in $B_{\frac{R}{2}}$. Furthermore, Proposition \ref{FL2} tell us that $\|w_k\|_{C^{2s+\alpha}}$ is uniformly bounded in $B_{\frac{R}{4}}$. We may do the same for a sequence $R_j\rightarrow \infty$, and pass to a diagonal subsequence (which we still denote as $\{w_k\}$) converging in $C_{\mbox{loc}}^{2s+\beta}({\mathbb{R}}^n)$, where $\beta\in(0,\alpha)$, to $\tilde{w}$ which satisfies
\begin{equation}\label{4.13}
\left\{\begin{array}{rl}
(-\Delta)^s{\tilde{w}}=0, \tilde{w}\leq 2^q ~~ \mbox{ in }{{\mathbb{R}}^n},\hspace{0.6cm}\\
\tilde{w}(0)=1.\hspace{3.8cm}
\end{array}\right.
\end{equation}
Since $\tilde{w}$ is $s$-harmonic and bounded above in ${{\mathbb{R}}^n}$, we have $\tilde{w}\equiv1$ in ${{\mathbb{R}}^n}$ by the Liouville theorem(\cite{8}). Therefore, $w_k\rightarrow 1$ in $C_{\mbox{loc}}({{\mathbb{R}}^n})$, which contradicts the second of (\ref{4.12}).
$\hfill\square$

Recalling that (\ref{asa}) plays an important role in our proof of Lemma \ref{FP2}. To ensure that (\ref{asa}) holds, the growth condition $u(x)=\mbox{O}(|x|^\tau)$ at infinity with $\tau<2s$ is necessary even if it seems to be very restricted.

Let $v$ be a solution of (\ref{1}) and we set
\begin{equation*}
\psi(x) =\frac{\Gamma(\frac{n-2s}{2})}{2^{2s}\pi^{\frac{n}{2}}\Gamma(s)}\int_{\mathbb{R}^n}|x-y|^{2s-n}v_+^\gamma(y)dy,~~~x\in\mathbb{R}^n.
\end{equation*}
It is obvious that $\psi\geq 0$. Next we can derive the representation formula of $v$ via $\psi$.
\begin{lem}\label{FP3}
For any nonconstant solution $v=v(x)\in \mathcal{L}_{2s}\cap C_{\mbox{loc}}^{2s+\alpha}(\mathbb{R}^n)$ of (\ref{1}), suppose that  $|v(x)|= \mbox{O}(|x|^\tau)$ at infinity with $\tau<2s$, then there exists a constant $\sigma>0$ such that
\begin{equation}\label{4.14}
v(x)=\psi(x)-\sigma.
\end{equation}
Moreover, the support of $v_+$ is compact and
\begin{equation}\label{s1}
\lim\limits_{|x|\rightarrow \infty}(v(x)+\sigma)|x|^{n-2s}=\frac{\Gamma(\frac{n-2s}{2})}{2^{2s}\pi^{\frac{n}{2}}\Gamma(s)}\int_{\mathbb{R}^n}v_+^\gamma(y)dy.
\end{equation}
\end{lem}
{\bf Proof } Firstly, we  show that $\psi=\psi(x)$ is well defined and
\begin{equation}\label{4.15}
\psi(x)\rightarrow 0 ~~~~\mbox{as}~ |x|\rightarrow + \infty.
\end{equation}
Lemma \ref{FP2} and the integral constraint in (\ref{1}) assure
\begin{equation}\label{4.16}
v_+\in L^\mu(\mathbb{R}^n)~~\mbox{for any} ~\mu\in[\frac{n(\gamma-1)}{2s},\infty].
\end{equation}
Given $R>0$, we introduce
$$ \psi_1(x)=\int_{|y-x|\geq R}|x-y|^{2s-n}v_+^\gamma(y)dy,$$
$$ \psi_2(x)=\int_{|y-x|< R}|x-y|^{2s-n}v_+^\gamma(y)dy,$$
that is $\psi=\frac{\Gamma(\frac{n-2s}{2})}{2^{2s}\pi^{\frac{n}{2}}\Gamma(s)}(\psi_1+\psi_2)$. At first, $\psi_1$ is estimated by
\begin{eqnarray}
\nonumber
&&{ \psi_1(x)\leq}
\begin{cases}
R^{2s-n}{\int_{|z|\geq R}}{v_+^\gamma}(x-z)dz,    &\gamma \in (1,\frac{n}{n-2s}],\\
\left({\int_{|z|\geq R}}|z|^{-n(1+\frac{2s}{(n-2s)\gamma-n})}\right)^{\frac{(n-2s)\gamma-n}{n(\gamma-1)}}\left({\int_{|z|\geq R}}|v_+|^{\frac{n(\gamma-1)}{2s}}(x-z)\right)^{\frac{2\gamma s}{n(\gamma-1)}},
 &\gamma\in(\frac{n}{n-2s},\frac{n+2s}{n-2s}),
\end{cases}
\\ &&\label{4.18}
~~~~~~~~\quad\indent\leq
\begin{cases}
R^{2s-n}\|v_+\|_\gamma^\gamma,& \gamma\in(1,\frac{n}{n-2s}],\\
R^{-\frac{2s}{\gamma-1}}C_1(n,\gamma)\|v_+\|_{\frac{n(\gamma-1)}{2s}}^\gamma,   &\gamma\in(\frac{n}{n-2s},\frac{n+2s}{n-2s}).
\end{cases}
\end{eqnarray}
Because of $\gamma \frac{n-s}{s} \geq \frac{n(\gamma-1)}{2s}$, it holds by (\ref{4.16}) that
\begin{align} \nonumber \label{4.17}
0\leq \psi_2(x)&=\int_{|y-x|< R}|x-y|^{2s-n}v_+^\gamma(y)dy\\  \nonumber
&\leq\left(\int_{|y-x|< R}|x-y|^{(2s-n) \frac{n-s}{n-2s}}dy\right)^{\frac{n-2s}{n-s}}\left(\int_{|y-x|< R}v_+^{\gamma \frac{n-s}{s}}(y)dy\right)^{\frac{s}{n-s}}\hspace{1.2cm}\\  \nonumber
&=\left(\int_{|y-x|< R}|x-y|^{s-n}dy\right)^{\frac{n-2s}{n-s}}\left(\int_{|y-x|< R}v_+^{\gamma \frac{n-s}{s}}(y)dy\right)^{\frac{s}{n-s}}\hspace{2.2cm}  \nonumber\\
&\leq C_1(n,R) \|v_+\|_{L^{\gamma \frac{n-s}{s}}(B_R(x))}^\gamma \rightarrow 0 ~~ \mbox{as}~   |x|\rightarrow +\infty.
\end{align}

Note that $\gamma\geq\frac{n(\gamma-1)}{2s}$ if $\gamma\in (1,\frac{n}{n-2s}]$. This  and (\ref{4.16}) yield that $\|v_+\|_\gamma<\infty$. From this and (\ref{4.16})-(\ref{4.17}), we deduce that $\psi$ is well defined. From (\ref{4.18}), we know that for any $\varepsilon>0$, there exists $R=R(\varepsilon)>0$ such that $0\leq \psi_1(x)\leq \frac{\varepsilon}{2}$ for any  $x\in\mathbb{R}^n$. And from (\ref{4.17}), we know that there exists $M>0$ such that $\psi_2(x)\leq \frac{\varepsilon}{2}$ for $|x|\geq M$. So (\ref{4.15}) follows.

\par Plainly we have
$$
(-\Delta)^s(v-\psi)=0~~ \mbox{in} ~{\mathbb{R}^n}.
$$
By Lemma \ref{FP2} and (\ref{4.15}) we have
$$\sup_{\mathbb{R}^n}(v-\psi)<+\infty. $$
Then, the Liouville theorem guarantees that there exists a real number $\sigma$ such that $v-\psi=-\sigma$. We claim that $\sigma>0$. Indeed, if not, then $v$ is a nonconstant nonnegative solution of $(-\Delta)^sv=v^\gamma$ in ${\mathbb{R}^n}$, which contradicts with the nonexistence of \cite{chen} in the subcritical case  $\gamma\in(1,\frac{n+2s}{n-2s})$.  Hence, we obtain (\ref{4.14}). Furthermore, (\ref{4.14}), (\ref{4.15}) and $\sigma>0$ yield that the support of $v_+$ is compact. This and the same argument as lemma \ref{L1999912} imply that (\ref{s1}) holds.
$\hfill\square$

\par Now, it is suffice to prove our main Theorem \ref{F1}.

{\bf Proof of Theorem \ref{F1}.}
Comparing (\ref{012}) with (\ref{4.14}), by replacing $p$ with $s$, we may apply almost the same argument of Theorem \ref{T1} to finish the proof of Theorem \ref{F1}, and we omit it.
$\hfill\square$

\section{Radial symmetry of solutions to fractional  equations without integral constraint}

{\bf Proof of Theorem \ref{F2}.}
\rm Without loss of generality, we assume that $u$ is not a constant function. For $x=(x_1,x')$, denote $\Sigma_\lambda=\{x|x_1<\lambda\}$, $T_\lambda=\{x|x_1=\lambda\}$, $x^\lambda=(2\lambda-x_1,x')$ and $w_\lambda(x)=u_\lambda(x)-u(x)$, we can derive that
\begin{equation}\label{x33}
(-\triangle)^s (u_\lambda(x)-u(x))+u_\lambda(x)-u(x)=g(u_\lambda)-g(u).
\end{equation}
Combined this with (\ref{1.6}), it is easy to verify that
\begin{equation*}
(-\triangle)^s w_\lambda(x)+w_\lambda(x)-Cu^\mu(x) w_\lambda(x)\geq0, ~~~~x\in \Sigma_\lambda^-,
\end{equation*}
where $\Sigma_\lambda^-=\{x\in\Sigma_\lambda|w_\lambda(x)<0\}$.
Therefore, at points where $w_\lambda$ is negative, we can obtain that
\begin{equation}\label{5.1}
(-\triangle)^s w_\lambda(x)+c(x) w_\lambda\geq0, ~~~x\in\Sigma_\lambda^-,
\end{equation}
where $c(x)=1-Cu^\mu(x)$. Clearly $c(x)$ is bounded from below.
\par $\mathbf{Step~ 1}$: We first apply Proposition $3.7^{(')}$  (Decay at infinity) to show that for sufficiently negative $\lambda$,
\begin{equation}\label{5.2}
w_\lambda(x)\geq0,~~~~ x\in\Sigma_\lambda.
\end{equation}
By condition (\ref{1.7}), we obtain that
\begin{equation*}
c(x)>-\tilde{C}|x|^{-2s},
\end{equation*}
and hence assumption (\ref{3.2}) in Proposition $3.7^{(')}$ is satisfied. Consequently, there exists $R_0>0$ such that for $\lambda<-R_0$,
\begin{equation*}
w_\lambda(x)\geq0,~~~~ x\in\Sigma_\lambda.
\end{equation*}
This verifies (\ref{5.2}).
\par $\mathbf{Step~ 2}$: Step 1 provides a starting point, from which we can now move the plane ${\mathrm{T}}_\lambda$ to the right as long as (\ref{5.2}) holds to its limiting position. Let
$$ \lambda_0=\sup\{\lambda|w_\mu(x)\geq0, \mbox{ for any } x\in\Sigma_\mu, \mu\leq\lambda\}.$$
Step 1 tells us that $\lambda_0\geq-R_0$.
We further claim that
$$\lambda_0\leq R_0.$$
Otherwise, if $\lambda_0>R_0$, then for any $\lambda\in(R_0, \lambda_0]$, from Proposition $3.7^{(')}$ and the rotation invariance of (\ref{x33}) we derive
\begin{equation}\label{www}
w_\lambda(x) \leq 0,~~ x \in \Sigma_\lambda.
\end{equation}
On the other hand, from the definition of $\lambda_0$, for  $\lambda\leq\lambda_0$ we have
\begin{equation}\label{eeee}
w_\lambda(x) \geq 0,~~ x \in \Sigma_\lambda.
\end{equation}
From (\ref{www})-(\ref{eeee}), we conclude that $w_\lambda(x)\equiv 0$ for $\lambda\in(R_0, \lambda_0]$. Hence $u$ is constant, which gives a contradiction.

\par $\mathbf{Step~ 3}$: We will show that $w_{\lambda_0}(x)\equiv0$ for $x\in\Sigma_{\lambda_0}$. Otherwise, the maximum principle for anti-symmetric functions (Proposition \ref{FL5}) yields that
\begin{equation}\label{5.5}
w_{\lambda_0}(x)>0, ~~~~ x\in\Sigma_{\lambda_0}.
\end{equation}

Again by Proposition $3.7^{(')}$, we obtain for any $\lambda\in(\lambda_0,\lambda_0+1)$, there exists $R_1>0$ such that
\begin{equation}\label{ooo}
w_{\lambda}(x)\geq 0,  ~~~x\in \Sigma_{\lambda}\setminus B_{R_1}.
\end{equation}
It follows from (\ref{5.5}) that there exists constant $\delta>0$ and $c_0>0$ such that
$$ w_{\lambda_0}(x)\geq c_0,~~~ x\in \overline{\Sigma_{\lambda_0-\delta}\cap B_{R_1}(0)}.$$
Since $w_\lambda$ depends on $\lambda$ continuously, there exists $\varepsilon>0$ and $\varepsilon< \delta$ such that for all $\lambda\in(\lambda_0,\lambda_0+\varepsilon)$, we have
\begin{equation}\label{conti}
w_\lambda(x)\geq 0,~~~ x\in \overline{\Sigma_{\lambda_0-\delta}\cap B_{R_1}(0)}.
\end{equation}

Since $c(x)$ is bounded from below, taking the narrow region principle (Proposition \ref{FL6}), we derive that there exists a small positive constant $\varepsilon(<\delta)$ such that $w_\lambda(x)\geq0$ in $(\Sigma_\lambda\backslash \Sigma_{\lambda_0-\delta})\cap B_{R_1}(0)$ for $\lambda\in(\lambda_0,\lambda_0+\varepsilon)$.
This and (\ref{ooo})-(\ref{conti}) yield that $w_\lambda(x)\geq0$  in $\Sigma_\lambda$  for $\lambda\in(\lambda_0,\lambda_0+\varepsilon)$, which  contradicts with the definition of $\lambda_0$. Therefore we must have $w_{\lambda_0}\equiv0$.
\par Since the problem is invariant with respect to rotation, we can take any direction as the $x_1$ direction. Hence we  have that $u$ is radially symmetric about some $x_0\in {\mathbb{R}}^n$.
$\hfill\square$

{\bf Proof of Theorem \ref{F3}.}
\rm Let $T_\lambda$, $x^\lambda, u_\lambda$, $w_\lambda(x)$, $\Sigma_\lambda$ and $\Sigma_\lambda^-$ be defined as in the previous section again.  Then
$$ (-\triangle)^s(u_\lambda(x)-u(x))=h(u_\lambda(x))-h(u(x)).$$
From (\ref{1.9}), we can obtain that
$$\left|\frac{h(u_\lambda(x))-h(u(x))}{u_\lambda(x)-u(x)}\right|\leq C (u_\lambda(x)+u(x))^\mu ,~~\mbox{for all } u_\lambda(x)\neq u(x),$$
then at points where $w_\lambda$ is negative, we have
$$ (-\triangle)^sw_\lambda(x)-C 2^\mu \cdot u^\mu(x)\cdot w_\lambda(x)\geq0.  $$

\par $\mathbf{Step~ 1}$: Set $c(x)=-C 2^\mu \cdot u^\mu(x)$, then it is obvious that $c(x)$ is bounded from below. Combined this with $u(x)=\mbox{O}(\frac{1}{|x|^m})$, we obtain that
$$ c(x) \sim -C 2^\mu \cdot \frac{1}{|x|^{m\mu}}.$$
And hence assumption (\ref{3.2}) in Proposition $3.7^{(')}$ is satisfied. Then there exists $R_0>0$ such that for $\lambda<-R_0$, we must have
\begin{equation}\label{5.6}
w_\lambda(x)\geq0,~~ x\in\Sigma_\lambda.
\end{equation}
\par $\mathbf{Step~ 2}$: Let
$$ \lambda_0=\sup\{\lambda|w_\mu(x)\geq0, \mbox{ for any } x\in\Sigma_\mu, \mu\leq\lambda\}.$$
We obtain that $\lambda_0\in[-R_0,R_0]$ by using the same method as Theorem \ref{F2}. We can also show that $w_{\lambda_0}(x)\equiv0$ for $x\in\Sigma_{\lambda_0}$.
In fact, if
\begin{equation}\label{5.8}
w_{\lambda_0}\geq0  \mbox{ and }  w_{\lambda_0}\not\equiv0, ~~x\in \Sigma_{\lambda_0},
\end{equation}
we must have
\begin{equation}\label{5.9}
w_{\lambda_0}(x)>0,~~x\in\Sigma_{\lambda_0},
\end{equation}
due to Proposition \ref{FL5}. Under our assumptions, we have
$$ \lim_{|x|\rightarrow \infty}w_\lambda(x)=0   \mbox{ and } c(x) \mbox{ is bounded from below }.$$
Combining this with the narrow region principle and the decay at infinity, through a similar argument as in the previous section, we derive there exists some $\varepsilon>0$ such that for any $\lambda\in(\lambda_0,\lambda_0+\varepsilon)$, $w_\lambda(x)\geq0$ for $x\in\Sigma_\lambda$, which is a contradiction.
\par Therefore, $w_{\lambda_0}(x)\equiv0$ for $x\in\Sigma_{\lambda_0}$. The results of this theorem follow.
$\hfill\square$

\end{document}